# HOMOGENEOUS ASYMPTOTIC LIMITS OF HAAR MEASURES OF SEMISIMPLE LINEAR GROUPS AND THEIR LATTICES

FRANÇOIS MAUCOURANT

ABSTRACT. We show that Haar measures of connected semisimple groups, embedded via a representation into a matrix space, have a homogeneous asymptotic limit when viewed from far away and appropriately rescaled. This is still true if the Haar measure of the semisimple group is replaced by the Haar measure of a irreducible lattice of the group, and the asymptotic measure is the same. In the case of an almost simple group of rank greater than 2, a remainder term is also obtained. This extends and precises anterior results of Duke, Rudnick and Sarnak, and Eskin-McMullen in the case of a group variety.

1. INTRODUCTION

1.1. **Statement of the results.** This paper deals with the asymptotic properties of the Haar measure of a semisimple Lie group, seen as a measure on a vector space $\text{End}(V)$ of endomorphisms of a vector space $V$. By a Radon measure, we mean a regular measure that is finite on compact sets. We prove here the following general result, which is a stronger form of the computation of the growth rate of a group that was done previously in some particular cases.

**Theorem 1.** *Let $G$ be a connected semisimple noncompact real Lie group with finite center, $\rho : G \to GL(V)$ a representation of $G$ with finite kernel in a finite dimensional real vector space $V$, and $\mu$ a Haar measure on $G$. Then there is a rational positive number $d$, an integer $e$ between $0$ and $\text{rank}_{\mathbf{R}} G - 1$, and a nonzero Radon measure $\mu_\infty$ on the space $\text{End}(V)$ of endomorphisms of $V$ such that we have : for any $f$ continuous of compact support from $\text{End}(V)$ to $\mathbf{C}$,*

$$\lim_{T \to +\infty} \frac{1}{T^d \ln(T)^e} \int_G f\left(\frac{\rho(g)}{T}\right) d\mu(g) = \int_{\text{End}(V)} f d\mu_\infty.$$

The numbers $d$ and $e$ admit the following geometrical description : if $\mathcal{C}$ is the convex hull of the weights of the representation $\rho$, then $d$ is the unique positive real number such that the half-sum of positive roots with multiplicities divided by $d/2$ lies on the boundary of $\mathcal{C}$, and $e$ the codimension minus one of the minimal face of the polyhedron $\mathcal{C}$ containing this point (see section 2.1). Moreover, the measure $\mu_\infty$ is *homogeneous of degree $d$* in the following sense : for any Borel set $E \subset \text{End}(V)$, any $t > 0$, we have $\mu_\infty(tE) = t^d \mu(E)$. It is of full support in a manifold $M$ which is an homogeneous space of $G \times G$, if the first copy of $G$ acts by matrix multiplication on the left and the second one on the right; moreover, the closure of $M$ in $\text{End}(V)$ might be smaller than the set of limits of sequences $(\rho(g_i)/T_i)_i$ (see the example of the Adjoint representation of $SL(3, \mathbf{R})$ in section 4.2). The function $T^d \ln(T)^e$ is called *the growth rate of $G$ relative to $\rho$*; there is a natural ordering on such functions given by the lexicographical order on $(d, e)$. It is not hard to show that the growth rate of a normal subgroup $H$ of $G$ with respect to $\rho_{|H}$ is always smaller or equal to the growth rate of the group $G$ with respect to $\rho$.

In [DRS], Duke, Rudnick and Sarnak investigated the properties of integer points of affine homogeneous varieties of simple noncompact algebraic groups; roughly speaking, they showed, among other things and under suitable assumptions, that the asymptotic of the number of integer points of norm less that $T$ is the same that the Haar measure of a ball of radius $T$; their proof was then simplified by Eskin and McMullen [EM] using a dynamical approach, and latter Eskin, Mozes and Shah [EMS] generalized this result to a wider setting than that of affine homogeneous varieties. Here we prove the following.

**Theorem 2.** *Let $G$, $\rho$ as in Theorem 1. Let $\Gamma$ be a lattice in $G$, and assume either that $\Gamma$ is irreducible, either that every nontrivial normal subgroup $H$ of $G$ has strictly smaller growth rate relative to $\rho_{|H}$ than $G$ with respect to $\rho$. Then for the same $d, e, \mu_\infty$ as in Theorem 1, for any $f$ continuous of compact support from $\text{End}(V)$ to $\mathbf{C}$, we have*

$$\lim_{T \to +\infty} \frac{1}{T^d \ln(T)^e} \sum_{\gamma \in \Gamma} f\left(\frac{\rho(\gamma)}{T}\right) = \frac{1}{Covol(\Gamma)} \int_{\text{End}(V)} f d\mu_\infty.$$

The condition on the subgroup growth have a geometrical description in terms of the dual of $\mathcal{C}$ (See Lemma 2.2).

If one drops both irreducibility and strictly smaller growth of proper normal subgroups, the conclusion of Theorem 2 can fail, as in the following example. Let $\Gamma = SL(2, \mathbf{Z}) \times SL(2, \mathbf{Z}) \subset G = SL(2, \mathbf{R}) \times SL(2, \mathbf{R})$, and the representation $\rho$ of $G$ on $V = \mathbf{R}^2 \otimes \mathfrak{sl}(2, \mathbf{R})$ acting on pure tensors by $\rho(g_1, g_2)v \otimes w = \rho_1(g_1)v \otimes \rho_2(g_2)w$, where $\rho_1$ is the standard representation and $\rho_2$ the Adjoint representation; here the growth rate $T^2$ of $G$ is the same as the growth rate of the first copy of $SL(2, \mathbf{R})$. In this case, there is a limit measure for the lattice, but it is different from the one for the group; see section 4.4 for details. It is different from the counter-examples with the same type of behaviour given in [DRS],[EM] because here the variety $\rho(G)$ is affine symmetric, whereas theirs was not.

The following corollary of Theorem 1 was announced independently by Gorodnik and Weiss in [GW], Theorem 2.7.

**Corollary 1.1.** *(Same assumptions as Theorem 1.) For any norm $||.||$ on $End(V)$, let $C$ be the $\mu_\infty$ measure of the ball of radius 1 around 0 for the chosen norm.*

Then, as $T$ tends to infinity,

$$\int_{g\in G\,,\,||\rho(g)||\leq T} dg \sim CT^d \ln(T)^e.$$

Theorem 2 implies the following Corollary.

**Corollary 1.2.** *(Same assumptions as Theorem 2.) For any norm $||.||$ on $End(V)$, let $C$ be the $\mu_\infty$ measure of the ball of radius 1 around 0 for the chosen norm. Then, as $T$ tends to infinity,*

$$\#\{\gamma \in \Gamma\,,\,||\rho(\gamma)|| \leq T\} \sim \frac{C}{covol(\Gamma)} T^d \ln(T)^e.$$

Duke, Rudnick and Sarnak showed in fact that Corollary 1.2 is a consequence of the Corollary 1.1 under the hypothesis that $G$ is almost simple. If $e = 0$, it is possible to improve Theorem 1 to give a remainder term.

**Theorem 3.** *(Same hypotheses and notations as Theorem 1.) We assume $e = 0$. There exists an effective $\alpha_0 > 0$ such that for any $\alpha \in ]0,1]$, for any Hölder map $f$ of exponent $\alpha$ and compact support, we have*

$$\frac{1}{T^d}\int_G f\left(\frac{\rho(g)}{T}\right) d\mu(g) = \int_{End(V)} f d\mu_\infty + O(T^{-\alpha\alpha_0}).$$

*Moreover, for any norm $||.||$ on $End(V)$, we have as $T$ tends to infinity,*

$$\int_{g\in G\,,\,||\rho(g)||\leq T} dg = CT^d + O(T^{d-\alpha_0}),$$

*where $C$ is the $\mu_\infty$ measure of the ball of radius 1 around 0 for the chosen norm.*

It is also possible to find a remainder term for Theorem 2 in the case when $G$ is almost simple of real rank greater than 2, following a suggestion of Babillot ([Bab], p47-48).

**Theorem 4.** *Same hypotheses and notations as Theorem 2. Assume $e = 0$, and that $G$ is almost simple of rank greater than 2. There exists a effective $\alpha_0 > 0$, independent of $\Gamma$, such that for any $\alpha \in ]0,1]$, for any Hölder map $f$ of exponent $\alpha$ and compact support, we have*

$$\sum_{\gamma\in\Gamma} f\left(\frac{\rho(\gamma)}{T}\right) = \frac{T^d}{covol(\Gamma)} \int f d\mu_\infty + O(T^{d-\alpha\alpha_0}).$$

*Moreover, for any norm $||.||$ on $End(V)$, we have as $T$ tends to infinity,*

$$\#\{\gamma \in \Gamma\,,\,||\rho(\gamma)|| \leq T\} = \frac{C}{covol(\Gamma)} T^d + O(T^{d-\alpha_0}),$$

*where $C$ is the $\mu_\infty$ measure of the ball of radius 1 around 0 for the chosen norm.*

This gives a partial answer to a question of Eskin-McMullen in the case of a group variety. It was already obtained by Duke, Rudnick and Sarnak [DRS] in the case of the standard representation of $SL(n, \mathbf{R})$; however in this particular case their remainder term $O(T^{n^2-n-1/(n+1)+\epsilon})$ is still better than the $O(T^{n^2-n-1/6n+\epsilon})$ for $n \geq 7$ that is computed here (see Proposition 4.1). Another known case of remainder term, due to Röttger ([R1],[R2]), and not covered by this Theorem, is the representation of $SL_2(\mathbf{O}_K)$ of $2 \times 2$ matrices with entries in the integer ring of a number field $K$ on $(\mathbf{R}^2)^r \times (\mathbf{C}^2)^s$.

If $e > 0$, then one obtains remainder terms of the form $O(T^d \ln T^{e-1})$. We give an example of application of Theorem 4 in 4.3, see Proposition 4.2.

Let us now describe the easiest example. If one takes $G = SL(2, \mathbf{R})$ imbedded in $\text{End}(\mathbf{R}^2)$ via the standard representation, then $d = 2$ and $e = 0$ as is well-known. The limit measure $\mu_\infty$ is supported on the manifold of rank 1 matrices which has the following parametrization, for $(\theta, \phi) \in [0, 2\pi[^2$ and $r > 0$:

$$\begin{bmatrix} r \sin \theta \sin \phi & r \sin \theta \cos \phi \\ r \cos \theta \sin \phi & r \cos \theta \cos \phi \end{bmatrix}.$$

In these coordinates, the limit measure $\mu_\infty$ can be written (up to a constant depending on the choice of Haar measure) $r dr \wedge d\theta \wedge d\phi$. The following corollary was originally the motivation for the study of homogeneous limits of Haar measures; compare with Ledrappier [L] (see also Nogueira [N] and Gorodnik [G]).

**Corollary 1.3.** *Let $f$ be a continuous map from $\mathbf{R}^2$ to $\mathbf{R}$. Let $\Gamma$ be a lattice in $SL(2, \mathbf{R}) \subset M(2, \mathbf{R})$, and $||.||_{M(2,\mathbf{R})}$ be the matrix norm $\sqrt{a^2 + b^2 + c^2 + d^2}$. We put*

$$\Gamma_T = \{\gamma \in \Gamma \; : \; ||\gamma||_{M(2,\mathbf{R})} \leq T\}.$$

*Define $D$ to be the disk of radius 1 and center 0 in $\mathbf{R}^2$, endowed with the usual euclidean norm. Then for any $v$ in $\mathbf{R}^2 - \{0\}$, we have*

$$\lim_{T \to +\infty} \frac{1}{\#\Gamma_T} \sum_{\gamma \in \Gamma_T} f(\gamma v/T) = \frac{2}{\pi^2} \int_{||v||D} f(y) \frac{\sqrt{||v||^2 - ||y||^2}}{||v||^2.||y||} dy.$$

*Proof.* (sketch) The limit given is nothing else than the (normalized) image measure of the limit measure restricted to the unit ball in $M(2, \mathbf{R})$, $1_{B(0,1)}\mu_\infty$, via the evaluation map from $M(2, \mathbf{R})$ to $\mathbf{R}^2 : M \mapsto M(v)$. □

E. Peyre pointed out to me the striking similarities between these results and a conjecture of Manin (see [FMT], [P]), where the number of points of bounded height on certain algebraic varieties has also a growth of the type $T^\alpha \ln(T)^\beta$; the exponents $\alpha, \beta$ are obtained by a very similar geometric construction.

1.2. **Organization of the paper.** In 2.1, we develop the notations and collect some general geometrical properties. In 2.2, we state the technical Propositions 2.1 and 2.2, of which Theorem 1 and 3, and Theorem 4 respectively, are easy consequences. Sections 2.3, 2.4, 2.5 are devoted to the proof of Proposition 2.1.

It is based on the formula for Haar measure using the Cartan decomposition $G = K \times \exp(\mathfrak{a}^+) \times K$, but we begin by ignoring the compact group at first, which does not change much the asymptotic (see 2.5). The main term is given by Laplace-type integrals defined on certain polyhedra (2.3). In section 2.6 we collect some analytical informations that will be used in the proof of Theorem 2 and Proposition 2.2. The proof of Theorem 2 that is given in 3.1 is a simplification of the argument of Duke-Rudnick-Sarnak [DRS] to our case, and rely, as in [EM], on the vanishing of matrix coefficients Theorem of Howe-Moore [HM]. Then, in 3.2, we state and prove a effective version of a Theorem due to Katok and Spatzier [KS] about decay of matrix coefficients for simple Lie groups of rank greater than 2, using the estimate for $K$-finite vectors due to Oh [O]. It is then used to prove along sections 3.3, 3.5 the Proposition 2.2, whose proof is actually an effective version of that of Theorem 2. We then give some examples, $SL(n, \mathbf{R})$ with the standard representation in 4.1, the Adjoint representation of $SL(3, \mathbf{R})$ in 4.2, the representation of $SL(3, \mathbf{R})$ on $\mathbf{R}^3 \oplus (\mathbf{R}^3)^*$ in 4.3, and lastly in 4.4 a counter-example to Theorem 2 when the lattice is reducible and does not satisfy the smaller normal subgroup growth hypothesis.

## 2. Asymptotics of smooth Haar Measures

**2.1. Notations and geometrical properties.** Let $G$ be a semisimple, noncompact, connected real Lie group with finite center. We will write $\mathfrak{g} = Lie(G)$ for its Lie algebra, and we fix $\mathfrak{g} = \mathfrak{k} + \mathfrak{p}$ a Cartan decomposition, where $\mathfrak{k} = Lie(K)$ is the Lie algebra of a maximal compact subgroup $K$. Let $\mathfrak{a} \subset \mathfrak{p}$ a Cartan subalgebra, $A$ the abelian, connected subgroup such that $Lie(A) = \mathfrak{a}$, and $r = dim(\mathfrak{a})$ the $\mathbf{R}$-rank of $G$. We will write $\Sigma$ for the nontrivial weights of the Adjoint representation, $(m_\alpha)_{\alpha \in \Sigma}$ their multiplicity, $W$ the Weyl group, $\Sigma^+$ a set of positive roots, $\langle .|. \rangle$ the duality brackets between $\mathfrak{a}^*$ and $\mathfrak{a}$, and
$$\mathfrak{a}^+ = \{a \in \mathfrak{a} \,|\, \forall \alpha \in \Sigma^+, \, \langle \alpha | a \rangle \geq 0\},$$
the positive, closed Weyl chamber associated to $\Sigma^+$. We put
$$\beta = \sum_{\alpha \in \Sigma^+} m_\alpha \alpha,$$
the sum of positive roots with multiplicities.

Given a faithful representation $\rho : G \to GL(V)$, where $V$ is a real vector space of dimension $n$, we write $\Phi = (\chi_1, .., \chi_n) \subset \mathfrak{a}^*$ for the weights of the representation $\rho$, counted with multiplicities. Let $\mathcal{C}$ the convex hull of the set $\Phi \subset \mathfrak{a}^*$; it is a convex and compact polyhedron, and we call
$$\mathcal{C}^* = \{a \in \mathfrak{a} \,|\, \forall \chi \in \mathcal{C}, \, \langle \chi | a \rangle \leq 1\},$$
the dual convex polyhedron in $\mathfrak{a}$. Recall that a *face* of a convex polyhedron $P$ is, either $P$ itself, either a nonempty intersection of $P$ with a supporting hyperplane of $P$. The *dimension* of a face is the dimension of the affine subspace it generates.

By (relative) *interior of a face* $\mathcal{F}$, we mean $\mathcal{F}$ minus the union of proper subfaces of $\mathcal{F}$. The polyhedron $\mathcal{C}^*$ is compact, thanks to the following Lemma.

**Lemma 2.1.** *The zero vector $0 \in \mathfrak{a}^*$ is in the interior of $\mathcal{C}$.*

*Proof.* Since $G$ is a unimodular group,

$$\sum_{\chi \in \Phi} \chi = 0, \tag{1}$$

thus 0 is a barycenter of elements of $\Phi$, so belongs to $\mathcal{C}$. If 0 belongs to one of the faces of dimension $(r-1)$ of $\mathcal{C}$, there is an $a \in \mathfrak{a} - \{0\}$ such that $\langle \chi | a \rangle \geq 0$ for all $\chi$ in $\mathcal{C}$. From Equation (1), this means that $\langle \chi | a \rangle = 0$ for all $\chi$ in $\mathcal{C}$, and since the kernel of $\rho$ is finite, $Vect(\Phi) = \mathfrak{a}^*$ so we have $a = 0$, which is a contradiction. □

Thus, we also have the following description of $\mathcal{C}$.

$$\mathcal{C} = \{\chi \in \mathfrak{a}^* \,|\, \forall\, a \in \mathcal{C}^*,\, \langle \chi | a \rangle \leq 1\},$$

Let $d > 0$ the smallest positive real number such that $\beta/d$ belongs to $\mathcal{C}$, it is well defined thanks to the preceding Lemma. Let $\mathcal{F}_\beta$ be the face of maximal codimension containing $\beta/d$. Let

$$e = \operatorname{codim}(\mathcal{F}_\beta) - 1.$$

Denote $\Delta = \mathfrak{a}^+ \cap \mathcal{C}^*$, and

$$\mathcal{F}_\beta^* = \{a \in \mathcal{C}^* \,|\, \langle \beta | a \rangle = d\}.$$

Recall (see for example [Be]) that the dual face of a $k$-dimensional face $\mathcal{F}$ is the $r - k - 1$-dimensional face

$$\mathcal{F}^* = \{a \in \mathcal{C}^* \,|\, \forall\, \chi \in \mathcal{F},\, \langle \chi | a \rangle = 1\},$$

but in fact one can prove that if $\chi$ is in the interior of the face $\mathcal{F}$, we also have $\mathcal{F}^* = \{a \in \mathcal{C}^* \,|\, \langle \chi | a \rangle = 1\}$. Thus $\mathcal{F}_\beta^*$ is the dual face of $\mathcal{F}_\beta$, and $e$ is then also the dimension of $\mathcal{F}_\beta^*$. Also, $\mathcal{F}_\beta^*$ is the set of point of $\mathcal{C}^*$ where $\beta$ attains its maximum, $d$ being the maximum of $\beta$ restricted to $\mathcal{C}^*$.

**Lemma 2.2.** *There is no proper ideal $\mathfrak{h}$ of $\mathfrak{g}$ containg $\mathcal{F}_\beta^*$ if and only if every nontrivial normal subgroup of $G$ has striclty smaller growth rate relative to $\rho$.*

*Proof.* Let $H$ be the connected normal subgroup with Lie algebra $\mathfrak{h}$ which is a ideal containg $\mathcal{F}_\beta^*$. Let $\mathfrak{a}', (\mathcal{C}')^*, \beta', d', e'$ be the same objects than $\mathfrak{a}, \mathcal{C}^*, \beta, d, e$, but with respect to $H$ rather than $G$. We can assume that $\mathfrak{a}' \subset \mathfrak{a}$, denote by $i$ this inclusion map, and then we have a canonical map $i^* : \mathfrak{a}^* \to (\mathfrak{a}')^*$. It is easily checked that $i^*(\beta) = \beta'$, and $(\mathcal{C}')^* = \mathcal{C}^* \cap \mathfrak{a}'$; since $\mathcal{F}_\beta^* \subset \mathfrak{a}'$, the definitions of $(d, e)$ in terms of maximum of $\beta$ on $\mathcal{C}^*$, and the dimension of the set where the maximum is attained, implies that $d = d'$ and $e = e'$. The reciprocal implication is similar. □

The following geometrical properties will be needed in the sequel.

**Lemma 2.3.** *The face $\mathcal{F}^*_\beta$ is contained in the closed positive Weyl chamber $\mathfrak{a}^+$.*

*Proof.* Let $a \in \mathcal{F}^*_\beta$, we have $\langle \beta | a \rangle = d$. Let $\alpha \in \Sigma^+$ a simple root, and $w_\alpha$ the associated involution in the Weyl group, which permutes simple roots other than $\alpha$ and such that $w_\alpha(\alpha) = -\alpha$. So $\beta - w_\alpha(\beta) = u\alpha$, where $u$ is a strictly positive integer. We have
$$\langle w_\alpha(\beta/d) | a \rangle = 1 - u\langle \alpha | a \rangle.$$
Since $\mathcal{C}$ is invariant under the Weyl group $W$, and $a$ belongs to $\mathcal{C}^*$, we obtain $\langle \alpha | a \rangle \geq 0$. Since this is true for all simple roots $\alpha$, this shows that $a$ belongs to $\mathfrak{a}^+$. □

For two linear forms $\gamma_1, \gamma_2$ in $\mathfrak{a}^*$, we will write $\gamma_1 \leq \gamma_2$ if this is true on restriction to $\mathfrak{a}^+$.

**Lemma 2.4.** *Let $\gamma$ an element of $\mathfrak{a}^*$ such that $w(\gamma) \leq \beta$ for all $w$ in the Weyl group $W$. Then $\gamma/d$ belongs to $\mathcal{C}$.*

*Proof.* We have to prove that for all $a \in \mathcal{C}^*$, $\langle \gamma/d | a \rangle \leq 1$. Let $a$ be in $\mathcal{C}^*$, there is a element $w$ in the Weyl group $W$ such that $wa$ is in $\mathfrak{a}^+$. So we have
$$\langle \gamma | a \rangle = \langle w(\gamma) | w(a) \rangle \leq \langle \beta | wa \rangle \leq d.$$
□

The proof of the following easy geometrical lemma is omitted.

**Lemma 2.5.** *Assume $\gamma$ lies in the interior of a face $\mathcal{F}$ of $\mathcal{C}$ of codimension at least 1, and let $\gamma'$ be in $\mathcal{C} - \mathcal{F}$. Then there is a $\lambda$ in $]0,1]$ such that $\lambda \gamma' + (1-\lambda)\gamma$ lies in the interior of a face of codimension stricly smaller than the codimension of $\mathcal{F}$.*

2.2. **General Estimates.** We now state the general technical estimates that implies Theorems 1, 3 and 4.

Let $\mathcal{B}$ be the space of measurable, bounded maps of compact support from the space $\text{End}(V)$ of endomorphisms of $V$ to $\mathbf{C}$, and $\mathcal{B}_1$ the subspace of $\mathcal{B}$ of maps which vanish outside the ball of center $0 \in \text{End}(V)$ of radius 1 for some norm on $\text{End}(V)$. If $f \in \mathcal{B}$, put $R_f$ for the least positive real number $r$ such that $f$ vanishes outside the ball of center $0$ and radius $r$; obviously, $f \in \mathcal{B}_1$ if and only if $R_f \leq 1$. We write $||f||_\infty$ for the supremum of $|f|$. In order to control the regularity of such maps, we introduce for $\delta > 0$ and $\epsilon \geq 0$ the following set
$$D_f(\delta, \epsilon) = \{g \in \text{End}(V) \mid \exists h, ||h - g|| \leq \delta, |f(g) - f(h)| > \epsilon\}.$$
This set will be of use to give estimates that make sense for both continuous functions and characteristic functions of set.

**Proposition 2.1.** *Same hypotheses than Theorem 1. There exists $d, e, \mu_\infty$, with the properties announced in Theorem 1, and there exists $\lambda \geq 0, \xi > 0, \xi' \geq 0, \tau \geq 0$, and a Radon measure $\nu \geq \mu_\infty$ homogeneous of degree $d$, and a constant $c > 1$* such that for any $f$ in $\mathcal{B}$, any $\delta$ in $]0,1[$, any $\epsilon$ in $[0,1[$ and any $T \geq c\delta^{-\lambda}/R_f$, we have

$$\left| \frac{1}{T^d \ln(T)^e} \int_G f\left(\frac{\rho(g)}{T}\right) d\mu(g) - \int_{\text{End}(V)} f d\mu_\infty \right|$$

$$\leq cR_f^d \left[ ||f||_\infty \left( R_f^{-d} \nu(D_f(cR_f\delta, \epsilon)) + \delta^{-\tau}(R_f T)^{-\xi} + \frac{(R_f T)^{-\xi'}}{\ln(TR_f)} + h\left(\delta, TR_f/c\right) \right) + \epsilon \right],$$

where $h(\delta, T) = 0$ if $e = 0$, and $h(\delta, T) = -\frac{\ln \delta + 1}{\ln T}$ otherwise. Moreover $\xi' > 0$ if $e = 0$.

The preceding Proposition implies Theorems 1 and 3 in the following way. We only sketch the calculations. For Theorem 1, take $\epsilon = \omega_f(c\delta)$ a continuity module for $f$, that is a map such that $D_f(\delta, \omega_f(\delta)) = \emptyset$ and such that $\omega_f(\delta)$ tends to zero as $\delta$ tends to zero; then let $T$ tends to infinity, and then $\delta$ to zero. This proves Theorem 1. In Theorem 3, we assumed $e = 0$, and for a Hölder map $f$ of exponent $\alpha$, we have $\omega_f(\delta) = \delta^\alpha$. Thus, we put $\epsilon = c^\alpha \delta^\alpha$; define

$$v = \inf(1/\lambda, \xi/(\alpha + \tau)),$$

and then put $\delta = T^{-v}$. It can be checked that $\tau v - \xi < 0$, and that $T \geq \delta^{-\lambda}$ as required, and then the result is given by

$$\alpha_0 = \inf(\xi - \tau v, \xi', v) > 0.$$

For the characteristic function $f$ of a ball $B(1)$ of radius 1 around 0 for some norm, one can take $\epsilon = 0$, then the measure of $D_f(c\delta, 0)$ is the measure of some small annulus $B(1 + c\delta) - B(1 - c\delta)$ with respect to some homogeneous measure $\nu$, and thus is less than some constant times $\delta$; thus the upper bound can be minimized with the constant given for Hölder maps of exponents $\alpha = 1$.

**Proposition 2.2.** *Let $G$ be a almost simple real Lie group of rank greater than 2, $\rho$ a representation on a finite dimensional vector space $V$, and let $\Gamma$ be a lattice in $G$. Let $d, e, \mu_\infty$ be given by Theorem 1. Let $\lambda \geq 0, \xi > 0, \xi' \geq 0, \tau \geq 0, \nu$ be given by Proposition 2.1. There exists $\xi_1 > 0, \tau_1 \geq 0$, and a constant $c > 0$ such that for any $f$ in $\mathcal{B}$, any $\delta$ in $]0,1[$, any $\epsilon$ in $[0,1[$ and any $T \geq c\delta^{-\lambda} R_f^\lambda$, we have*

$$\left| \frac{1}{T^d \ln(T)^e} \sum_{\gamma \in \Gamma} f\left(\frac{\rho(\gamma)}{T}\right) - \frac{1}{covol(\Gamma)} \int_{\text{End}(V)} f d\mu_\infty \right|$$

$$\leq cR_f^d [||f||_\infty (R_f^{-d} \nu(D_f(cR_f\delta, \epsilon)) + \delta^{-\tau}(R_f T)^{-\xi} + \delta^{-\tau_1}(R_f T)^{-\xi_1}$$

$$+ \frac{(R_f T)^{-\xi'}}{\ln(TR_f)} + h\left(\delta, TR_f/c\right)) + \epsilon],$$

*where $h(\delta, T) = 0$ if $e = 0$, and $h(\delta, T) = -\frac{\ln \delta + 1}{\ln T}$ otherwise. Moreover $\xi' > 0$ if $e = 0$.*

This Proposition implies Theorem 4 in exactly the same way Proposition 2.1 implies Theorem 3.

2.3. **Asymptotics of some measures defined by linear forms.** We fix a basis of $V$ in which $\exp(\mathfrak{a})$ is diagonal, that is

$$\rho \circ \exp(a) = Diag(\exp\langle\chi_1|a\rangle, .., \exp\langle\chi_n|a\rangle).$$

The vector space $\mathrm{End}(V)$ is endowed with the operator norm associated to the supremum norm in $V$ (with respect to the chosen basis). For $f$ in $\mathcal{B}$, we will note $f_d$ the map $f$ restricted to the set of diagonal matrices, that is

$$f_d(h_1, .., h_n) = f(Diag(h_1, .., h_n)).$$

We denote by $da$ the choice of a nonzero Lebesgue measure on $\mathfrak{a}$. If $\mathcal{F}$ is a proper face of $\mathcal{C}$, we define $\mathcal{M}_\mathcal{F}$ to be the image of $\mathfrak{a}$ by the map $a \mapsto (f_i(a))_{1 \leq i \leq n}$, where $f_i(a) = \exp\langle\chi_i|a\rangle$ if $\chi_i \in \mathcal{F}$, 0 otherwise. It is a manifold of dimension $\dim\mathcal{F} + 1$.

Before we state the main result of this section, let us state a easy upper bound.

**Lemma 2.6.** *Let $f$ in $\mathcal{B}$ and $\gamma$ in $\mathfrak{a}^*$. Let $u$ be the maximum of $\gamma$ restricted to $\Delta = \mathcal{C}^* \cap \mathfrak{a}^+$. Then for all $T \geq 1$, we have*

$$\left|\int_{\mathfrak{a}^+} f\left(\frac{\rho \circ \exp(a)}{T}\right) \exp\langle\gamma|a\rangle da\right| \leq \mathcal{L}(\Delta)||f||_\infty (R_f)^u T^u (\ln(TR_f))^r,$$

*where $\mathcal{L}(\Delta)$ stands for the Lebesgue measure of $\Delta$ : $\mathcal{L}(\Delta) = \int_\Delta da$.*

*Proof.* First, since $f(\frac{\rho\circ\exp(a)}{T}) = f_d(\exp\langle\chi_1|a\rangle/T, .., \exp\langle\chi_n|a\rangle/T)$, this quantity vanishes when $a$ in $\mathfrak{a}^+$ is outside $\ln(TR_f)\Delta$. So we have

$$\left|\int_{\mathfrak{a}^+} f\left(\frac{\rho \circ \exp(a)}{T}\right) \exp\langle\gamma|a\rangle da\right| \leq ||f||_\infty \int_{\ln(TR_f)\Delta} \exp\langle\gamma|a\rangle da.$$

Then apply a change of variable $b = a/\ln(TR_f)$ to obtain the expected result. $\square$

The rest of the section will be devoted to the proof of the following Proposition.

**Proposition 2.3.** *Let $k \geq 1$, and $\gamma_1, .., \gamma_k$ be distinct, nonzero linear forms in $\mathfrak{a}^*$, such that*

- *for all $i$, $\gamma_i \leq \gamma_1$,*
- *there exists a positive real number $\zeta$, and a face $\mathcal{F}$ of $\mathcal{C}$ of codimension at least one, containing $\gamma_i/\zeta$ for all $i$,*
- *the element $\gamma_1/\zeta$ lies in the interior of $\mathcal{F}$, which means that $\mathcal{F}$ is the face of $\mathcal{C}$ of minimal dimension containing $\gamma_1/\zeta$.*

*Let $s$ be the dimension of $\mathcal{F}$, and $t = (t_1, .., t_k)$ a nonzero vector in $\mathbf{R}^k$. There there is a linear form $L$ on $\mathcal{B}$, a Radon measure $L_+$ (both may be zero), and constants $C, \lambda \geq 0, \tau \geq 0, \xi > 0$ such that for any $f$ in $\mathcal{B}$, any $\delta, \epsilon$ in $]0,1[$ et $[0,1[$ respectively, and any $T > \delta^{-\lambda}/R_f$, we have*

$$\left|\frac{1}{T^\zeta \ln(T)^{r-1-s}} \int_{\mathfrak{a}^+} f\left(\frac{\rho \circ \exp(a)}{T}\right) \sum_{i=1}^k t_i \exp\langle\gamma_i|a\rangle da - L(f)\right|$$

$$\leq \left[ C||f||_\infty \left( R_f^{-\zeta} L_+(D_f(R_f\delta, \epsilon)) + \delta^{-\tau}(TR_f)^{-\xi} + h(\delta, TR_f) \right) + C\epsilon \right] R_f^\zeta,$$

where $h(\delta, T) = 0$, if $s + 1 = r$, and otherwise $h(\delta, T) = \frac{-\ln(\delta)+1}{\ln(T)}$. *The linear form $L$ is a linear combination of measures in the Lebesgue class of $\mathcal{M}_\mathcal{F}$, and $L_+$ is also in the Lebesgue class of this manifold. Moreover, $L$ is nonzero if and only if we have the following inclusion $\{a \in \mathcal{C}^* \mid \langle \gamma_1 | a \rangle = \zeta\} \subset \mathfrak{a}^+$ (equivalently $\mathcal{F}^* \subset \mathfrak{a}^+$). In the case $L$ is nonzero, then $L$ and $L^+$ are both homogeneous of degree $\zeta$.*

Taking $\delta = 1/2$, $\epsilon = 0$, we have the following.

**Corollary 2.1.** *Same assumptions and notations as Proposition 2.3. There is a constant $c > 0$ such that for any $f$ in $\mathcal{B}$, any $T \geq 2/R_f$, we have*

$$\left| \int_{\mathfrak{a}^+} f\left(\frac{\rho \circ \exp(a)}{T}\right) \sum_{i=1}^k t_i \exp\langle \gamma_i | a \rangle da \right| \leq c ||f||_\infty (TR_f)^\zeta (\ln(TR_f))^{r-1-s}.$$

□

First, we can assume that $R_f > 0$, and if we put $f'(x) = f(R_f x)$, then $f'$ belongs to $\mathcal{B}_1$, and it is easily checked that

$$D_{f'}(\delta, \epsilon) = \frac{1}{R_f} D_f(R_f \delta, \epsilon).$$

If $s+1 = r$, using this reduction, it is easily checked that is is sufficient to prove the Proposition when $R_f = 1$; in the case $r > s+1$, it is a little more complicated because one have to use the Corollary when $R_f = 1$ and the particular form of $h$; we skip the details. Thus it is enough to prove Proposition 2.3 for $f \in \mathcal{B}_1$, and from now on we will make this assumption.

Since $f$ vanishes outside the ball of center 0 and radius 1, we have

$$I(T) = \int_{\mathfrak{a}^+} f\left(\frac{\rho \circ \exp(a)}{T}\right) \sum_{j=1}^k t_j \exp\langle \gamma_j | a \rangle da = \int_{\ln(T)\Delta} f\left(\frac{\rho \circ \exp(a)}{T}\right) \sum_{j=1}^k t_j \exp\langle \gamma_j | a \rangle da.$$

Let $\mathcal{F}'$ a codimension 1 face containing $\mathcal{F}$. We can find a set of $s+1$ weights, say $\chi_1, .., \chi_{s+1}$ up to changing the indices, which are a affine basis of the affine subspace of dimension $s$ containing $\mathcal{F}$. Thus, there exists real numbers $\mu_{i,j}$ (which are allowed to be negative), such that for all $j = 1, .., k$, we have $\sum_{i=1}^{s+1} \mu_{i,j} = 1$ and

$$(2) \qquad \gamma_j/\zeta = \sum_{i=1}^{s+1} \mu_{i,j} \chi_i.$$

Since $\mathcal{F}'$ is of codimension 1 and does not contain 0, we can freely change the indices of weights other than $\chi_1, .., \chi_{s+1}$ such that $(\chi_1, .., \chi_r)$ is a basis of $\mathfrak{a}^*$, contained in $\mathcal{F}'$. Let $(\chi_1^*,..,\chi_r^*)$ be the dual basis of $\mathfrak{a}$. The following lemma will be useful in the sequel.

**Lemma 2.7.** *The element $\sum_{i=1}^r \chi_i^*$ belongs to $\mathcal{C}^*$.*

*Proof.* The dual face $(\mathcal{F}')^*$ of $\mathcal{F}'$ is a 0-dimensional face (a point), and using the definition this must be $\sum_{i=1}^r \chi_i^*$, so is in $\mathcal{C}^*$. $\square$

Let $D > 0$ be the real number such that
$$da = D d\chi_1 \wedge ... \wedge d\chi_r.$$

For $i = 1,..,s+1$, put

$$(3) \qquad x_i(a) = \frac{\exp\langle \chi_i | a \rangle}{T},$$

and for $i = 1,..r-s-1$,

$$(4) \qquad y_i(a) = \frac{\langle \chi_{s+1+i} | a \rangle}{\ln(T)}.$$

Remark that if $s = r - 1$, the variables $y_i$ are not defined and should be skipped in all the sequel. Let $E_T$ be the image of $\ln(T)\Delta$ by the map
$$a \mapsto (x_1(a),..x_{s+1}(a), y_1(a),..,y_{r-s-1}(a)).$$

For $i = 1,..s+1$, we can differentiate to obtain $d\chi_i = dx_i/x_i$, and for $i = 1,..,r-s-1$, we have $d\chi_{s+1+i} = \ln(T)dy_i$. Therefore, the change of variable $a \mapsto (x_i, y_i)$ writes

$$(5) \qquad \frac{I(T)}{T^\zeta(\ln(T))^{r-s-1}} = D \int_{E_T} f_d\left(P_1(\bar{x})T^{l_1(\bar{y})},..,P_n(\bar{x})T^{l_n(\bar{y})}\right) \sum_{j=1}^k t_j \prod_{i=1}^{s+1} x_i^{\zeta\mu_{i,j}-1} d\bar{x}d\bar{y},$$

where $\bar{x}, \bar{y}$ denotes the vectorial variables $(x_1,..,x_{s+1})$ and $(y_1,..y_{r-s-1})$, $P_i$ are the expressions

$$(6) \qquad P_i(\bar{x}) = \prod_{j=1}^{s+1} x_j^{\langle \chi_i | \chi_j^* \rangle},$$

and $l_i$ the affine forms

$$(7) \qquad l_i(\bar{y}) = \left(\sum_{j=1}^{s+1} \langle \chi_i | \chi_j^* \rangle\right) - 1 + \sum_{j=1}^{r-s-1} \langle \chi_i | \chi_{j+s+1}^* \rangle y_i.$$

Let us describe the integration set $E_T$. Translating the two conditions $\chi_i \leq \ln(T)$ and $\alpha \geq 0$ defining $\ln(T)\Delta$, we get that the set $E_T \subset (\mathbf{R}_+^*)^{s+1} \times \mathbf{R}^{r-s-1}$ is defined by the following inequalities, for $i = 1,..,n$

$$\text{(8)} \qquad T^{l_i(\bar{y})} P_i(\bar{x}) \leq 1,$$

and for all $\alpha$ in $\Sigma^+$,

$$\text{(9)} \qquad T^{l_\alpha(\bar{y})} P_\alpha(\bar{x}) \geq 1,$$

with the convention that

$$P_\alpha(\bar{x}) = \prod_{j=1}^{s+1} x_j^{\langle \alpha | \chi_j^* \rangle},$$

and

$$l_\alpha(\bar{y}) = \sum_{j=1}^{s+1} \langle \alpha | \chi_j^* \rangle + \sum_{j=1}^{r-s-1} \langle \alpha | \chi_{j+s+1}^* \rangle y_j.$$

We define

$$Y_\infty = \{\bar{y} \in \mathbf{R}^{r-s-1} \mid \text{for } i = 1,..,n,\ l_i(\bar{y}) \leq 0,\ \text{for } \alpha \in \Sigma^+,\ l_\alpha(\bar{y}) \geq 0\},$$

$$X_\infty = \{\bar{x} \in [0,1]^{s+1} \mid \text{for } i = 1,..,n \text{ s.t. } l_i \equiv 0,\ P_i(\bar{x}) \leq 1,\ \text{for } \alpha \in \Sigma^+ \text{ s.t.} l_\alpha \equiv 0,\ P_\alpha(\bar{x}) \geq 1\}.$$

**Lemma 2.8.** *There is a compact set $[0,1]^{s+1} \times [-C,1]^{r-s-1}$ containing $E_T$ for all $T \geq 1$. Let $E_\infty$ be the set of couples $(\bar{x}, \bar{y})$ that belongs to $E_T$ for all $T$ sufficiently large (it can happen that $E_\infty$ is empty). Then we have*

$$E_\infty = (X_\infty \times Y_\infty) \cap (]0,1]^{s+1} \times \mathbf{R}^{r-s-1}).$$

*The sets $Y_\infty$ and $X_\infty$ are compact, and $X_\infty$ is not of empty interior in $[0,1]^{s+1}$.*

*Proof.* The equality $E_\infty = (X_\infty \times Y_\infty) \cap (]0,1]^{s+1} \times \mathbf{R}^{r-s-1})$. is completly elementary, from Equations (8) and (9), and the fact that Equation (8) asserts for $i = 1,..,s+1$ that $x_i \leq 1$. From their definitions, $X_\infty$ and $Y_infty$ are closed, so $X_\infty$ is compact; moreover since $\Delta$ is compact, it follows from their definitions that the variables $y_i$ are bounded independently of $T$, so $Y_\infty$ is also compact. Equation (8) implies that if $(\bar{x}, \bar{y}) \in E_T$ then $\bar{x} \in X_\infty$; since $E_T$ is not of empty interior, neither is $X_\infty$. □

Let $M(\bar{x}) = (h_i)_{1 \leq i \leq n}$, where $h_i = P_i(\bar{x})$ if $l_i \equiv 0$, and $h_i = 0$ otherwise.

**Lemma 2.9.** *For all $\bar{x} \in \mathbf{R}_+^{s+1}$, and $\lambda > 0$ we have*

$$M(\lambda \bar{x}) = \lambda M(\bar{x}).$$

*Moreover, for almost all $(\bar{x}, \bar{y})$ in $E_\infty$, we have*

$$\lim_{T \to +\infty} (P_i(\bar{x}) T^{l_i(\bar{y})})_i = M(\bar{x}).$$

*Proof.* Remark that the set of $i$ such that $l_i$ is constant equal to zero is exactly the set such that $P_i$ is homogeneous of degree one, because of Equations (6) and (7), so $M(\lambda \bar{x}) = \lambda M(\bar{x})$. For $(\bar{x}, \bar{y})$ in $E_\infty$ and outside the affine hyperplanes defined by $l_i = 0$ for those $i$ such that $l_i$ is not identicaly zero, we have that $T^{l_i(\bar{y})} P_i(\bar{x})$ tends to zero because $l_i(\bar{y}) < 0$. □

Define
$$X_\omega = \{\bar{x} \in \mathbf{R}_+^{s+1} \,|\, \text{for } \alpha \in \Sigma^+ \text{ such that } l_\alpha \equiv 0, \ P_\alpha(\bar{x}) \geq 1\}.$$

Remark that the set of $\alpha$ such that $l_\alpha$ is constant equal to zero is exactly the set of $\alpha$ such that $P_\alpha$ is homogeneous of degree 0, thus $X_\omega$ is a cone, and in fact $X_\omega = \mathbf{R}_+ X_\infty$. Now we define for any $g$ in $\mathcal{B}$,

$$(10) \qquad L(g) = D\mathcal{L}(Y_\infty) \int_{X_\omega} g_d(M(\bar{x})) \sum_{j=1}^{k} t_j \prod_{i=1}^{s+1} x_i^{\zeta \mu_{i,j} - 1} d\bar{x}.$$

Here $\mathcal{L}(Y_\infty)$ stands for the Lebesgue measure of $Y_\infty \subset \mathbf{R}^{r-s-1}$, and is equal to 1 in the case $r = s+1$. It is not obvious that $L$ is well-defined, because some exponents $\mu_{i,j}$ might be nonpositive; however, note that the integrating set can be restricted to the compact set $R_g X_\infty$, because $g_d \circ M$ vanishes outside this set.

The problem of integration in the neighborhood of 0 will be the subject of lemma 2.11, but let us first state the following easy Lemma, whose proof is omitted.

**Lemma 2.10.** *Let $(v_1, .., v_k)$ be a finite set of vectors in a real vector space. Assume a vector $w$ lies in the interior of the convex hull of $(v_1, .., v_k)$. Then there exists $\kappa_1 > 0, .., \kappa_k > 0$ satisfying the following*

$$w = \sum_{i=1}^{k} \kappa_i v_i, \text{ and } \sum_{i=1}^{k} \kappa_i = 1.$$

□

**Lemma 2.11.** *There exists $\tau_1 > -1, .., \tau_{s+1} > -1$,*

$$Q(\bar{x}) = \prod_{i=1}^{s+1} x_i^{\tau_i},$$

*such that for all $\bar{x}$ such that there exists $\bar{y}$ and $T$ with $(\bar{x}, \bar{y}) \in E_T$, we have for any $j = 1, .., k$,*

$$Q(\bar{x}) \geq \prod_{i=1}^{s+1} x_i^{\zeta \mu_{i,j} - 1}.$$

*and thus, the integral defining $L$ is convergent, because $X_\infty$ is a subset of $[0, 1]^{s+1}$ and $\tau_i > -1$.*

*Proof.* It is sufficient to prove this for $j = 1$, because the inequality for $a \in \mathfrak{a}^+$, $\exp\langle\gamma_j|a\rangle \leq \exp\langle\gamma_1|a\rangle$ rewrites, after the change of variables,

$$T^\zeta \prod_{i=1}^{s+1} x_i^{\zeta\mu_{i,j}} \leq T^\zeta \prod_{i=1}^{s+1} x_i^{\zeta\mu_{i,1}}.$$

Let $\mathcal{R} = \Phi \cap \mathcal{F}$ be the set of weights included in $\mathcal{F}$. Since $\gamma_1/\zeta$ is in the interior of $\mathcal{F}$, which is the convex hull of the set $\mathcal{R}$, application of Lemma 2.10 yields positive real numbers $(\kappa_\chi)_{\chi \in \mathcal{R}}$ such that $\sum \kappa_\chi = 1$ and

$$\gamma_1/\zeta = \sum_{\chi \in \mathcal{R}} \kappa_\chi \chi.$$

Thus, combining with Equation (2) and splitting this sum between those $\chi$ appearing in $(\chi_1, .., \chi_{s+1})$ and the others, we obtain

$$\sum_{j=1}^{s+1} (\mu_{j,1} - \kappa_{\chi_j})\chi_j = \sum_{\forall j,\, \chi \neq \chi_j} \kappa_\chi \chi.$$

Since for all $a$ in $\ln(T)\Delta$ we have $\langle\chi|a\rangle \leq \ln(T)$, we can write for all $a$ in $\ln(T)\Delta$,

$$\sum_{j=1}^{s+1} (\mu_{j,1} - \kappa_{\chi_j})\langle\chi_j|a\rangle \leq \left(\sum_{\forall j,\, \chi \neq \chi_j} \kappa_\chi\right) \ln(T),$$

so

$$\prod_{j=1}^{s+1} (Tx_j)^{(\mu_{j,1} - \kappa_{\chi_j})} \leq T^{\sum_{\forall j,\, \chi \neq \chi_j} \kappa_\chi}.$$

The variable $T$ disappears because $\sum \mu_i = \sum \kappa_\chi = 1$. Taking the inverse inequality, we obtain

$$\prod_{j=1}^{s+1} x_j^{(\kappa_{\chi_j} - \mu_{j,1})} \geq 1.$$

We can rewrite it

$$\prod_{j=1}^{s+1} x_j^{\zeta\kappa_{\chi_j} - 1} \geq \prod_{j=1}^{s+1} x_j^{\zeta\mu_{j,1} - 1}.$$

All the exponents of the left hand side product being strictly bigger than $-1$, this concludes the lemma. $\square$

Remark that if $g$ is in $B_1$, then

$$L(g) = D\mathcal{L}(Y_\infty) \int_{X_\infty} g_d(M(\bar{x})) \sum_{j=1}^k t_j \prod_{i=1}^{s+1} x_i^{\zeta\mu_{i,j} - 1} d\bar{x},$$

because $g_d$ is zero whenever $P_i(\bar{x}) \geq 1$ for some $i$ such that $l_i \equiv 0$.

**Lemma 2.12.** *The form $L$, if nonzero, is homogeneous of degree $\zeta$.*

*Proof.* This follows from Equation (10), Lemma 2.9, and the fact that for all $j$, $\prod_{i=1}^{s+1} x_i^{\zeta \mu_{i,j}-1}$ is homogeneous of degree $\zeta - (s+1)$. □

**Lemma 2.13.** *Assume that*
$$\{a \in \mathcal{C}^* \mid \langle \gamma | a \rangle = \zeta\} \subset \mathfrak{a}^+,$$
*then the form $L$ is nonzero.*

*Proof.* The linear combination of product of powers $\sum_{i=1}^{k} t_j \prod_{i=1}^{s+1} x_i^{\zeta \mu_{i,j}-1}$ never vanishes on a open set of $\mathbf{R}^{s+1}$ if the family of exponents $((\zeta \mu_{i,j} - 1)_i)_j$ are distincts, and the vector $(t_1, .., t_k)$ is nonzero. It it the case since the $\gamma_i$ are distincts. So, since $X_\omega$ is never of empty interior (nor $X_\infty$), it is sufficient to show that $Y_\infty$ is in this case, not of empty interior. For $\bar{y} = (y_1, .., y_{r-s-1}) \in \mathbf{R}^{r-s-1}$, define
$$a(\bar{y}) = \sum_{i=1}^{s+1} \chi_i^* + \sum_{j=1}^{r-s-1} y_j \chi_{s+1+j}^*.$$

Then for $i = 1, .., n$, we have $l_i(\bar{y}) = \langle \chi_i | a(\bar{y}) \rangle - 1$, and for all $\alpha$ in $\Sigma^+$, $l_\alpha(\bar{y}) = \langle \alpha | a(\bar{y}) \rangle$. So $\bar{y}$ is in $Y_\infty$ if and only if $a(\bar{y})$ is in $\Delta$. Remark that the set $\{a \in \mathcal{C}^* \mid \langle \gamma_1 | a \rangle = \zeta\}$ is, by duality between $\mathcal{C}$ and $\mathcal{C}^*$, of dimension $r - s - 1$, and is exactly the set of $a$ in $\mathcal{C}$ such that $\langle \chi_i | a \rangle = 1$ for $i = 1, .., s+1$; so elements of this set can all be written $a(\bar{y})$ for some $\bar{y}$. Thus $Y_\infty$ is a polyedron of dimension $r - s - 1$, so has non-empty interior. □

Now, since all the properties of $L$ where checked, it remains to show the inequality expressed in Proposition 2.3. We can restrict ourselves to the case of $k = 1$, with only one $\gamma = \gamma_j$ for some $j$, and $t_1 = 1$, because the result will follow by linearity and the triangle inequality, if we put

$$L_+(g) = D\mathcal{L}(Y_\infty) \int_{X_\omega} g_d(M(\bar{x})) \sum_{j=1}^{k} |t_j| \prod_{i=1}^{s+1} x_i^{\zeta \mu_{i,j}-1} d\bar{x}.$$

Of course, $\gamma_j$ does not lie necessarily in the interior of $\mathcal{F}$, but we won't use this property again, it was only useful to show that $L$ was well defined. Remark also that, at this point, if one is only interested in showing that if $L(f)$ is nonzero and $f$ continuous, then the integral $I(T)$ is equivalent to $T^d (\ln(T))^{r-s-1} L(f)$, the result can be obtained here by using directly the dominated convergence Theorem to the expression (5), thanks to Lemmas 2.8 and 2.9.

We first fix the following constants. In the following formula, the sum runs through the set of index $i$ such that $l_i$ is a constant $l_i < 0$; in case this set in empty, take $\lambda = 0$.

$$\lambda = \sup_{i, l_i < 0,\ l_i \text{ const.}} -\frac{1}{l_i}. \tag{11}$$

We have to show that $\lambda \geq 0$. If $l_i$ is constant, from Equation (7) we have in particular

$$l_i = \left\langle \chi_i \Big| \sum_{j=1}^r \chi_j^* \right\rangle - 1,$$

and since $\sum_{j=1}^r \chi_j^* \in \mathcal{C}^*$ (Lemma 2.7), this implies that $l_i \leq 0$.

Write $\mathcal{T}$ for all couples $(i,j)$ with $i = 1, .., n$ such that $l_i$ is constant, $l_i \neq 0$, and $j \in \{1, .., s+1\}$ such that $\langle \chi_i | \chi_j^* \rangle < 0$.

$$\tau = \sup_{(i,j) \in \mathcal{T}} -\frac{\tau_j + 1}{(s+1)\langle \chi_i | \chi_j^* \rangle}, \quad \xi = \inf_{(i,j) \in \mathcal{T}} \frac{l_i(\tau_j + 1)}{(s+1)\langle \chi_i | \chi_j^* \rangle}. \tag{12}$$

The constant $\tau$ is nonnegative by definition of $\mathcal{T}$, and $\xi$ is positive for the same reason that $\lambda$. Note that these constants in fact depends on the choice of $\mathcal{F}'$ containig $\mathcal{F}$, of the elements of the basis $(\chi_i)$, and also on the choice of $\tau_i > -1$ in Lemma 2.11, so might be improved by wiser choices.

For now, we will try to give a upper bound to the difference

$$\mathcal{S}(T) = \frac{1}{D} \left| \frac{I(T)}{T^\zeta (\ln(T))^{r-s-1}} - L(f) \right|.$$

For some $\delta > 0$, write

$$E_{T,\delta} = \left\{ (\bar{x}, \bar{y}) \in E_\infty \mid \text{for all } i = 1, .., n,\ P_i(\bar{x}) T^{l_i(\bar{y})} < \delta,\ \text{or } l_i(\bar{y}) \equiv 0. \right\}.$$

Write $M_0(\bar{x}, \bar{y}) = diag(M(\bar{x})) \in \text{End}(V)$. We have the following inclusion.

$$E_T \cup E_\infty \subset (E_{T,\delta} - M_0^{-1}(D_f(\delta, \epsilon))) \cup (M_0^{-1}(D_f(\delta, \epsilon)) \cap E_\infty) \cup (E_\infty - E_{T,\delta}) \cup (E_T - E_\infty).$$

Thus, Equation (5) allows us to give the following bound, where $\epsilon \geq 0$:

$$\mathcal{S}(T) \leq \int_{E_{T,\delta} - M_0^{-1}(D_f(\delta,\epsilon))} \left| f_d\left((P_i(\bar{x})T^{l_i(\bar{y})})_{i=1,..,n}\right) - f_d(M(\bar{x})) \right| \prod_{i=1}^{s+1} x_i^{\zeta \mu_i - 1} d\bar{x} d\bar{y}$$

$$+ 2\|f\|_\infty \int_{M_0^{-1}(D_f(\delta,\epsilon)) \cap E_\infty} \prod_{i=1}^{s+1} x_i^{\zeta \mu_i - 1} d\bar{x} d\bar{y} + 2\|f\|_\infty \int_{(E_\infty - E_{T,\delta})} Q(\bar{x}) d\bar{x} d\bar{y}$$

$$+ \|f\|_\infty \int_{(E_T - E_\infty)} Q(\bar{x}) d\bar{x} d\bar{y}.$$

Let us call $\mathcal{S}_1(T,\delta,\epsilon), \mathcal{S}_2(T,\delta,\epsilon), \mathcal{S}_3(T,\delta), \mathcal{S}_4(T)$ the three terms of the sum on the right hand side. First, because of the definitions of $E_{T,\delta}, D_f(\delta,\epsilon)$ and $M(\bar{x})$, we have for all $(\bar{x},\bar{y}) \in E_{T,\delta} - M_0^{-1}(D_f(\delta,\epsilon))$ the following inequality

$$\left| f_d\left(P_1(\bar{x})T^{l_1(\bar{y})},..,P_n(\bar{x})T^{l_n(\bar{y})}\right) - f(M_0(\bar{x},\bar{y}))\right| \leq \epsilon.$$

Thus, since $E_{T,\delta} \subset E_\infty$, there is a constant $C_1 > 0$ such that

$$\mathcal{S}_1(T,\delta,\epsilon) \leq \epsilon \int_{E_\infty} \prod_{i=1}^{s+1} x_i^{\zeta\mu_i - 1} d\bar{x}d\bar{y} = C_1 \epsilon.$$

For the second term, we compute

$$\mathcal{S}_2(T,\delta,\epsilon) = 2||f||_\infty \int_{E_\infty} 1_{D_f(\delta,\epsilon)}(M_0(\bar{x},\bar{y})) \prod_{i=1}^{s+1} x_i^{\zeta\mu_i - 1} d\bar{x}d\bar{y} \leq 2||f||_\infty L(D_f(\delta,\epsilon))$$

For the two remaining terms, we will do separetly the two cases $r = s+1$ and $r > s+1$.

**First case :** $r = s+1$.
In this case, there is no $\bar{y}$ variables, so the affine forms $l_i$ are nonpositive constants. From the definition of $E_{T,\delta}$, we have

$$E_\infty - E_{T,\delta} \subset \left\{\bar{x} \in [0,1]^r \mid \text{for some } i \text{ such that } l_i \neq 0,\ T^{l_i} P_i(\bar{x}) > \delta\right\}.$$

For $j = 1,..,s+1$ and $\tau > 0$, we will write $G(j,\eta)$ for the set $\{\bar{x} \in [0,1]^r \mid x_j > \eta\}$. Let $\bar{x} \in E_\infty - E_{T,\delta}$, there is some $i$ with $l_i \neq 0$ such that

(13) $$T^{l_i} \prod_{j=1}^r x_j^{\langle \chi_i | \chi_j^* \rangle} \geq \delta.$$

Let $k$ such that $x_k^{\langle \chi_i | \chi_k^* \rangle}$ is maximum. Then

$$x_k^{r\langle \chi_i | \chi_k^* \rangle} \geq \prod_{j=1}^r x_j^{\langle \chi_i | \chi_j^* \rangle} \geq \delta T^{-l_i}.$$

Since $T \geq \delta^{-\lambda}$ by hypothesis, and $\delta < 1$, we have $T^{-l_i}\delta \geq \delta^{1-\lambda l_i} > 1$ because of Equation (11). Combining with $x_k \leq 1$, we conclude that $\langle \chi_i | \chi_k^* \rangle < 0$, that is $(i,k) \in \mathcal{T}$. Thus we have shown that :

$$E_\infty - E_{T,\delta} \subset \bigcup_{(i,j) \in \mathcal{T}} G\left(j, (\delta T^{-l_i})^{\frac{1}{r\langle \chi_i | \chi_j^* \rangle}}\right).$$

Now, a simple integration show that
$$\int_{G(j,\eta)} Q(\bar{x})d\bar{x} = \left(\prod_{i=1}^{r}(\tau_i+1)\right)^{-1}\eta^{\tau_j+1}.$$

Thus we obtain a first upper bound for $\mathcal{S}_3(T,\delta)$
$$\mathcal{S}_3(T,\delta) \leq \int_{E_\infty - E_{T,\delta}} Q(\bar{x})d\bar{x} \leq \left(\prod_{i=1}^{r}(\tau_i+1)\right)^{-1} \sum_{(i,j)\in\mathcal{T}} (\delta T^{-l_i})^{\frac{\tau_j+1}{(s+1)\langle\chi_i|\chi_j^*\rangle}},$$

Using the definitions of $\tau$ and $\xi$ (Equation (12)), we obtain the final upper bound
$$\mathcal{S}_3(T,\delta) \leq 2||f||_\infty rn \left(\prod_{i=1}^{r}(\tau_i+1)\right)^{-1} \delta^{-\tau} T^{-\xi}.$$

Here we have $E_T \subset E_\infty$. So in this case
$$\mathcal{S}_4(T) = 0.$$

**Second case :** $r > s+1$.

Let $(\bar{x},\bar{y})$ be in $E_\infty - E_{T,\delta}$. Then, for some $i$ with $l_i(\bar{y}) \neq 0$, we have
$$T^{l_i(\bar{y})} \prod_{j=1}^{s+1} x_j^{\langle\chi_i|\chi_j^*\rangle} \geq \delta.$$

Remark that $l_i(\bar{y}) < 0$ because otherwise, the left hand side would be greater than 1 for large $T$, contradicting the fact that $(\bar{x},\bar{y})$ is in $E_\infty$. Assume first that $l_i$ is a constant affine function. Because of the definition of $E_{T,\delta}$, this constant cannot be 0, and it is nonpositive. Thus, as we have seen, this type of inequality, together with the definition of Equation (11) implies that for some $j$ satisfying $\langle\chi_i|\chi_j^*\rangle < 0$, we have
$$x_j \leq (\delta T^{-l_i})^{\frac{1}{(s+1)\langle\chi_i|\chi_j^*\rangle}},$$
and, if we define $G(j,\eta)$ the set $\{\bar{x} \in [0,1]^{s+1} \mid x_j > \eta\}$, we have that
$$\bar{x} \in G(j,(\delta T^{-l_i})^{\frac{1}{(s+1)\langle\chi_i|\chi_j^*\rangle}}).$$

Now assume that $l_i$ is non-constant, from
$$1 \geq T^{l_i(\bar{y})} \prod_{j=1}^{s+1} x_j^{\langle\chi_i|\chi_j^*\rangle} \geq \delta,$$
we have
$$l_i(\bar{y}) \in \left[-\frac{\ln(P_i(\bar{x}))}{\ln T} + \frac{\ln\delta}{\ln T}, -\frac{\ln(P_i(\bar{x}))}{\ln T}\right],$$

which is an interval of width $\frac{\ln \delta}{\ln T}$, for fixed $\bar{x}$. Since $l_i$ is affine and non constant, there is a constant $C_i$, depending on $l_i$ and the compact set $Y_\infty$, such that for all real number $z$ and $\eta > 0$, the Lebesgue measure of the set $\{\bar{y} \in Y_\infty \mid l_i(\bar{y}) \in [z, z+\eta]\}$ is less that $C_i \eta$. We obtain the final upper bound

$$\mathcal{S}_3(T, \delta) \leq 2||f||_\infty \left(\prod_{i=1}^{s+1}(\tau_i + 1)\right)^{-1} \left((s+1)n\delta^{-\tau}T^{-\xi} + (\max_{i, l_i \text{ non-constant}} C_i)\frac{\ln \delta}{\ln T}\right).$$

Now, let $(\bar{x}, \bar{y}) \in E_T - E_\infty$. Thus, we have four possibilities : either for some $i$, $l_i(\bar{y}) > 0$, or for some $\alpha$, $l_\alpha(\bar{y}) < 0$, or $P_i(\bar{x}) > 1$ for some $i$ such that $l_i = 0$, either $P_\alpha(\bar{x}) < 1$ for some $\alpha$ such that $l_\alpha = 0$. As in the first case, the last two possibilities cannot occur. Let us treat the first case. Since $(\bar{x}, \bar{y}) \in E_T$, Equation (8) implies that

$$0 < l_i(\bar{y}) \leq -\frac{\ln P_i(\bar{x})}{\ln T}.$$

Applying the same argument as before, the set of $(\bar{x}, \bar{y})$ such that the preceding inequality occurs for this precise $i$ is less than

$$\int_{X_\infty} C_i \left(-\frac{\ln P_i(\bar{x})}{\ln T}\right) Q(\bar{x}) d\bar{x},$$

that is some constant times $1/\ln T$, because the integral is convergent thanks to Lemma 2.11. The second case is similar, we have by Equation (9),

$$0 > l_\alpha(\bar{y}) \geq -\frac{\ln P_\alpha(\bar{x})}{\ln T},$$

and this occurs for a set of $(\bar{x}, \bar{y})$ of measure at most some constant times $1/\ln T$. Thus, there is a constant $c > 0$ such that

$$\mathcal{S}_4(T) \leq \frac{c}{\ln T}.$$

This concludes the proof of Proposition 2.3.

2.4. **Sum of asymptotics.** We will note $h_i$ the nonzero integer coefficients and $\Omega$ the finite subset of $\mathfrak{a}^*$ such that

$$\prod_{\alpha \in \Sigma^+} (2\sinh\langle\alpha|a\rangle)^{m_\alpha} = \sum_{\gamma \in \Omega} h_\gamma \exp\langle\gamma|a\rangle.$$

**Lemma 2.14.** *For all $\gamma \in \Omega - \{0\}$, $\gamma/d$ belongs to $\mathcal{C}$.*

*Proof.* Observe that $\Omega$ is invariant under the Weyl group $W$, because the absolute value of the product $\prod_{\alpha \in \Sigma^+}(2\sinh\langle\alpha|a\rangle)^{m_\alpha}$ is invariant under the Weyl group. Since for any element $\gamma \in \Omega$, we have $\gamma \leq \beta$, Lemma 2.4 applies. □

Let $\iota_T : \text{End}(V) \to \text{End}(V)$ be the map $\iota_T(x) = x/T$.

**Proposition 2.4.** *Let $\nu$ be the measure defined for $f$ continuous of compact support on $End(V)$,*

$$\nu(f) = \int_{\mathfrak{a}^+} f(\rho \circ \exp(a)) \prod_{\alpha \in \Sigma^+} (2\sinh\langle \alpha | a \rangle)^{m_\alpha} da.$$

*Then there is nonzero Radon measures $\nu_\infty, L$, homogeneous of degree $d$, and constants $C > 0, \xi \geq 0, \lambda \geq 0$ and $\tau \geq 0$ such that for any $f$ in $\mathcal{B}$, any $\delta, \epsilon$ in $]0,1[$ et $[0,1[$ respectively, and any $T > \delta^{-\lambda}/R_f$, we have*

$$\left| \frac{\nu(f \circ \iota_T)}{T^d \ln(T)^e} - \nu_\infty(f) \right| \leq C R_f^d [||f||_\infty (R_f^{-d} L(D_f(R_f \delta, \epsilon))$$

$$+ \delta^{-\tau}(TR_f)^{-\xi} + (R_f T)^{-\xi'} \ln(T)^{-1} + h(\delta, TR_f)) + \epsilon],$$

*where $h(\delta, T) = 0$, if $s + 1 = r$, and otherwise $h(\delta, T) = \frac{-\ln(\delta)}{\ln(T)}$. Moreover, $\xi' > 0$ provided that $e = 0$.*

*Proof.* Let $\Omega_1$ be the set of $\gamma$ in $\Omega$ such that $\gamma/d$ belongs to the face $\mathcal{F}_\beta$, and $\Omega_2$ the complementary set in $\Omega$. If we define

$$S_i(f, T) = \int_{\mathfrak{a}^+} f\left( \frac{\rho \circ \exp(a)}{T} \right) \sum_{\gamma \in \Omega_i} h_\gamma \exp\langle \gamma | a \rangle da,$$

we can split the sum

$$\nu(f \circ \iota_T) = S_1(f, T) + S_2(f, T).$$

Proposition 2.3 applies directly to $S_1$ if one take $\zeta = d$, $\gamma_1 = \beta$, and Lemma 2.3 insures that the obtained limit linear form $\nu_\infty$ is nonzero. Let us prove that the second term is negligible, more precisely that there is a constant $C > 0$ such that

(14) $$S_2(f, T) \leq C||f||_\infty R_f^d T^{d-\xi'} \ln(T)^{e-1}.$$

For each $\gamma$ in $\Omega_2$, since $\gamma/d$ belongs to $\mathcal{C}$ thanks to Lemma 2.14, Lemma 2.5 applied to $\beta$ in $\mathcal{F}_\beta$ and $\gamma$, shows that there exists another linear form $\theta_\gamma$ such that $\gamma \leq \theta_\gamma$ and moreover $\theta_\gamma$ belongs to the interior of a face of codimension strictly smaller than $\mathcal{F}_\beta$. Thus Corollary 2.1 applied to $\gamma$ yields the required estimate (Equation 14). The fact that $\nu_\infty$ is a measure, i.e positivity, is a consequence of the positivity of $\iota_T^* \nu$. □

**2.5. Convolution with compact groups.** Here we give the final step of proof of Theorem 1. The Haar measure $\mu$ of $\rho(G)$ satisfies the following formula (see [H]), for a map $f$ in $\mathcal{B}$.

$$\int_{End(V)} f d\mu = \int_{K \times A^+ \times K} f(\rho(k \exp(a) k')) \prod_{\alpha \in \Sigma^+} (2\sinh\langle \alpha | a \rangle)^{m_\alpha} dk da dk'.$$

This can be written

$$\int_{\text{End}(V)} f d\mu = \int_{K \times K} \left( \int_{\text{End}(V)} f(kxk')d\nu(x) \right) dk dk',$$

and so we have

$$\mu(f \circ \iota_T) = \int_{K \times K} \left( \int_{\text{End}(V)} f \circ \iota_T(kxk')d\nu(x) \right) dk dk',$$

Let us define the measure $\mu_\infty$ by

$$(15) \qquad \int_{\text{End}(V)} f d\mu_\infty = \int_{K \times K} \left( \int_{\text{End}(V)} f(kxk')d\nu_\infty(x) \right) dk dk'.$$

Let $(k,k')f$ denotes the function $x \mapsto f(kxk')$. Since $K$ is compact, there is a constant $c > 0$ such that for all $k, k'$, and all $f$, we have

$$D_{(k,k')f}(\delta, \epsilon) \subset k^{-1} D_f(c\delta, \epsilon) k'^{-1},$$

and such that $R_{(k,k')f} \leq c R_f$. Thus, we have for all $\epsilon \in [0,1[$, $\delta \in ]0,1[$, and $T \geq c R_f \delta^\lambda$,

$$\left| \frac{\mu(f \circ \iota_T)}{T^d (\ln(T))^e} - \mu_\infty(f) \right| \leq C ||f||_\infty \left( \int_{K \times K} \nu_\infty(k^{-1} D_f(c\delta, \epsilon) k'^{-1}) dk dk' + h(\delta, T/(cR_f)) \right) + C\epsilon.$$

This concludes the proof of Proposition 2.1. Let us prove the remark stated in the introduction that $\mu_\infty$ could be seen as the unique (up to a multiplicative constant) $G \times G$-invariant measure on a certain $G \times G$-orbit. The measure $\nu_\infty$ is obtained (by construction) by integrating on the image of the map $M$, which is a $A$-orbit, and the boundary of this orbit has zero measure because of the integrability of the expression. One then constructs $\mu_\infty$ by integrating $\nu_\infty$ under the right and left $K$-actions, thus the measure $\mu_\infty$ is obtained by integrating on a $G \times G$-orbit an integrable expression, so the boundary of this orbit has measure zero. Since, on the other hand, $\mu_\infty$ is obtained as a limit of $G \times G$-invariant measure, it is also invariant under $G \times G$, but an orbit has at most one invariant measure class, so $\mu_\infty$ is necessarily the unique ergodic, invariant measure supported whose support is the closure of this orbit.

2.6. **Useful estimates.** For the proof of Theorem 2, we will need the following Proposition, whose technique of proof is very close to the one used in proof of Theorem 1. For a semisimple group $G$, we will say that $(g_i)_{i \geq 0}$ *tends strongly to infinity* if for all connected normal subgroups $H \neq G$ of $G$, the sequence $(g_i H)_i$ in $G/H$ tends to infinity.

**Proposition 2.5.** *Let $\xi$ be a bounded, non-negative map from $G$ to $\mathbf{R}$. Assume $\xi$ has the following property : for all sequence $(g_i)_i$ of $G$ tending strongly to infinity,* *then $\xi(g_i)$ tends to zero, and assume moreover that strict normal subgroups have strictly smaller growth. Then for all $f$ in $\mathcal{B}$,*

$$\int_G f(g/T)\xi(g)d\mu = o\left(T^d(\ln(T))^r\right).$$

*Proof.* (Sketch) First, we can assume that $\xi$ is $K$-invariant on the left and on the right, because

$$\xi'(g) = \sup_{(k,k')\in K^2} \xi(kgk'),$$

also tends to zero when $g$ tends strongly to infinity, because of the compactness of $K$. We have that $\exp(a)$ tends strongly to infinity, if and only if for all nontrivial ideals $\mathfrak{h}$, $a + \mathfrak{h}$ tends to infinity in $\mathfrak{a}/\mathfrak{h} \neq \mathfrak{a}$.

Using the $K$-invariance of $\xi$, we have to bound the following integral.

$$\int_{K\times\mathfrak{a}^+\times K} f(k\rho(\exp(a))k'/T)\xi(\exp(a)) \prod_{\alpha\in\Sigma^+} (2\sinh\langle\alpha|a\rangle)^{m_\alpha} dk\,da\,dk'.$$

Using the same convolution argument as in 2.5, and since the product of sinh is less that some multiple of $\exp\langle\beta|a\rangle$, the problem can be reduced to consider the following integral.

$$\int_{\mathfrak{a}^+} f(\rho\circ\exp(a)/T)\xi(\exp(a)) \prod_{\alpha\in\Sigma^+} \exp\langle\alpha|a\rangle da.$$

As usual, we can assume freely that $f$ is in $\mathcal{B}_1$. Applying the change of variables used in Proposition 2.3, we obtain a expression of the form, with $s = r - e - 1$,

$$T^d(\ln(T))^e \int_{E_T} f_d(P_1(\bar{x})T^{l_1(\bar{y})},...,P_n(\bar{x})T^{l_n(\bar{y})}) \prod_{i=1}^{s+1} x_i^{\zeta\mu_i-1}\xi(\exp(a(\bar{x},\bar{y},T)))d\bar{x}d\bar{y},$$

with

$$a(\bar{x},\bar{y},T) = \ln(T)\left(\sum_{i=1}^{s+1}\chi_i^* + \sum_{i=1}^{r-s-1} y_i\chi_{i+s+1}^*\right) + \sum_{i=1}^{s+1}\ln(x_i)\chi_i^*.$$

Since $\sum_{i=1}^{s+1}\chi_i^* + \sum_{i=1}^{r-s-1} y_i\chi_{i+s+1}^*$ describes $\mathcal{F}_\beta^*$ as $\bar{y}$ is in $Y_\infty$, Lemma 2.2 implies that $a(\bar{x},\bar{y},T)$ modulo any non-trivial ideal $\mathfrak{h}$ of $\mathfrak{g}$ goes to infinity as $T$ goes to infinity, for almost every $(\bar{x},\bar{y})$ in $X_\infty \times Y_\infty$. Thus for almost every $(\bar{x},\bar{y})$ in $X_\infty \times Y_\infty$, $\xi(\exp(a(\bar{x},\bar{y},T)))$ tends to zero as $T$ tends to infinity. Thus, applying the dominated convergence Theorem, the integral tends to zero as $T$ tends to infinity. □

**Lemma 2.15.** *If $g = kak'$ is a Cartan decomposition of $g$ in $G$, $k, k' \in K$ and $a \in A$, let $\xi(kak') = \exp(-\langle\gamma|a\rangle)$ for some linear form $\gamma$. Let $u$ be the maximum of $\beta - \gamma$ on restriction to $\Delta$. Then there exists a $c > 0$ such that for all $f$ in $\mathcal{B}$, and $T \geq 2/R_f$, we have*

$$\int_G f(g/T)\xi(g)d\mu \leq c||f||_\infty (R_fT)^u(\ln(TR_f))^r.$$

*Proof.* By the convolution argument, we are left with a integral of the form

$$\int_{\mathfrak{a}^+} f(\rho \circ \exp(a)/T) \prod_{\alpha \in \Sigma^+} \exp\langle \alpha - \gamma | a \rangle da,$$

and we apply Lemma 2.6. □

## 3. Lattices

3.1. **Proof of Theorem 2.** Let $G$ be a noncompact, semisimple, connected, Lie group of finite center, and $\rho : G \to GL(V)$ be a faithful representation on a finite dimensional real vector space. Let $\Gamma$ be a lattice in $G$; for convenience, we normalize the Haar measure of $G$ such that the covolume of $\Gamma$ is equal to 1. Here $\Gamma$ is assumed to be irreducible, or that strict normal subgroups have strictly smaller growth with respect to the representation $\rho$.

For a bounded, measurable map $f$ of compact support, from $G$ to $\mathbf{C}$, we associate the following map from $G \times G$ to $\mathbf{C}$.

$$\Psi(f)(g_1, g_2) = \sum_{\gamma \in \Gamma} f(g_1^{-1} \gamma g_2).$$

Since $\Gamma$ is discrete and $f$ of compact support, the sum is in fact finite for each $(g_1, g_2)$. It is clear that $\Psi(f)$ is left $\Gamma \times \Gamma$-invariant, thus can be seen as a map from $(\Gamma \backslash G)^2$ to $\mathbf{C}$; moreover it is easily checked that $\int_{(\Gamma \backslash G)^2} \Psi(f) = \int_G f$, and so $\Psi(f)$ is in $L^1((\Gamma \backslash G)^2)$. Thus, in order to prove Theorem 2, it is sufficient to prove that for all $f$ continuous of compact support from $\text{End}(V)$ to $\mathbf{C}$,

$$\lim_{T \to +\infty} \frac{\Psi(f \circ \iota_T)(\Gamma, \Gamma)}{T^d \ln(T)^e} = \int_{\text{End}(V)} f d\mu_\infty,$$

the map $f \circ \iota_T$ being seen here as a map from $G$ to $\mathbf{C}$ through the restriction to the image of $\rho$.

For a continuous map of compact support $\alpha$ from $(\Gamma \backslash G)^2$ to $\mathbf{C}$, we define for $h \in G$

$$\Psi^*(\alpha)(h) = \int_{\Gamma \backslash G} \alpha(\Gamma g, \Gamma g h) dg.$$

Since $\Gamma \backslash G$ is of finite volume, this is a bounded map.

**Proposition 3.1.** *The maps $\Psi$ and $\Psi^*$ are formally dual operators in the following sense. For bounded, measurable maps of compact support $f$ and $\alpha$ from $G$ to $\mathbf{C}$ and from $(\Gamma \backslash G)^2$ to $\mathbf{C}$ respectively, we have*

$$\int_{(\Gamma \backslash G)^2} \Psi(f)(\Gamma g_1, \Gamma g_2) \alpha(\Gamma g_1, \Gamma g_2) dg_1 dg_2 = \int_G f(h) \Psi^*(\alpha)(h) dh.$$

*Proof.* Let $D$ be a fundamental domain for $\Gamma$ in $G$. Then

$$\int_{D \times D} \Psi(f)(\Gamma g_1, \Gamma g_2) \alpha(\Gamma g_1, \Gamma g_2) dg_1 dg_2 = \int_{D \times G} f(g_1^{-1} g_2') \alpha(\Gamma g_1, \Gamma g_2') dg_1 dg_2',$$

where we have put the change of variable $g'_2 = \gamma g_2$, $g'_2 \in G$, $\gamma \in \Gamma$ and $g_2 \in D$. The latter expression is also equal to

$$\int_{(g_1,h) \in D \times G} f(h) \alpha(\Gamma g_1, \Gamma g_1 h) dg_1 dh,$$

if we use the change of variable $(g_1, h) = (g_1, g_1^{-1} g_2)$. □

For $1 > \eta > 0$, we fix a continuous map $\phi_\eta$ from $\Gamma \backslash G$ to $\mathbf{R}$, such that $\phi_\eta$ is nonnegative and zero outside a ball of radius $\eta$ around $\Gamma$ in $\Gamma \backslash G$ (the distance on $\Gamma \backslash G$ chosen being induced by a left-invariant Riemannian distance on $G$) and $\int_{\Gamma \backslash G} \phi_\eta = 1$. Now we put

$$\alpha_\eta(\Gamma g_1, \Gamma g_2) = \phi_\eta(\Gamma g_1) \phi_\eta(\Gamma g_2),$$

Let us recall the following Theorem.

**Theorem 5.** *(Howe-Moore [HM], [Z]) Let $G$ be a connected semisimple Lie group with finite center, $\pi$ a unitary representation of $G$ on a Hilbert space $\mathcal{H}$ such that if $H \neq \{1\}$ is a normal subgroup of $G$, then $\pi$ is without $H$-invariant vector. Then for any $v, w \in \mathcal{H}$,*

$$\lim_{g \to +\infty} \langle \pi(g) v, w \rangle = 0.$$

The following folklore corollary, whose proof is omitted, will be also useful.

**Corollary 3.1.** *Let $G$ be a connected semisimple Lie group with finite center, $\pi$ a unitary representation of $G$ on a Hilbert space $\mathcal{H}$ without $G$-invariant vector. Let $(g_i)_{i \geq 0}$ be a sequence in $G$ tending strongly to infinity (see 2.6 for definition). Then for any $v, w \in \mathcal{H}$,*

$$\lim_{i \to +\infty} \langle \pi(g_i) v, w \rangle = 0.$$

□

The next lemma is similar in spirit to Theorem 1.2 of [EM]. A basic observation is that $G$ acts unitarily and without fixed vectors on the Hilbert space $\mathcal{H}$ of $L^2$-functions on $\Gamma \backslash G$ of vanishing integral. Moreover, if $\Gamma$ is irreducible, there is no $H$-invariant vectors for $H$ a non-trivial normal subgroup.

**Lemma 3.1.** *For all $\eta > 0$, $|\Psi^* \alpha_\eta(h) - 1|$ tends to zero as $h$ tends strongly to infinity. If the lattice $\Gamma$ is irreducible, $|\Psi^* \alpha_\eta(h) - 1|$ tends to zero as $h$ tends to infinity.*

*Proof.*

$$|\Psi^* \alpha_\eta(h) - 1| = \left| \int_{\Gamma \backslash G} \phi_\eta(\Gamma g) \phi_\eta(\Gamma g h) dg - 1 \right| = |\langle \phi_\eta - 1, \pi(h^{-1})(\phi_\eta - 1) \rangle|,$$

and thus, Corollary 3.1 and Theorem 5 concludes each case. □

**Lemma 3.2.** *For any continuous $f$ of compact support, as $T$ tends to infinity,*
$$\lim_{T\to+\infty} \frac{1}{T^d(\ln(T))^e} \int_{(\Gamma\backslash G)^2} \Psi(f\circ\iota_T)\alpha_\eta = \int_{\mathrm{End}(V)} f d\mu_\infty.$$

*Proof.* We have
$$\int_{(\Gamma\backslash G)^2} \Psi(f\circ\iota_T)(g,g')\alpha_\eta(g,g')dgdg' = \int_G (f\circ\iota_T)(h)\Psi^*\alpha_\eta(h)dh$$
$$= \int_G (f\circ\iota_T)(h)dh + \int_G (f\circ\iota_T)(h)(\Psi^*\alpha_\eta(h)-1)dh.$$
The first integral is handled by Theorem 1. The second one is bounded by $o(T^d(\ln(T))^e)$, thanks to Lemma 3.1 and Proposition 2.5. □

Let $B(1_G,\eta)$ be the closed ball of radius $\eta$ in $\Gamma\backslash G$ for the chosen Riemannian metric on $G$.

In the compact neigborhood $B(1_G,1)$ of $1_G$ in $G$, the distance given by the operator norm on $\mathrm{End}(V)$ is equivalent to the Riemannian distance : there is some $c>0$ such that for all $(g,g')\in B(1_G,1)^2$,

(16) $\quad 1/c||\rho(g)-\rho(g')||_{\mathrm{End}(V)} \leq d_G(g,g') \leq c||\rho(g)-\rho(g')||_{\mathrm{End}(V)}.$

Let $f$ be continuous of support in the ball of center $0\in\mathrm{End}(V)$ and radius 1, and $\omega_f(\delta)$ a continuity modulus for $f$, that is a map such that $D_f(\delta,\omega_f(\delta))=\emptyset$, and $\lim_{\delta\to 0}\omega_f(\delta)=0$. We put

(17) $\quad\quad\quad\quad f_\eta^+(x) = \sup_{(g,g')\in B(1_G,\eta)^2} f(g^{-1}xg'),$

(18) $\quad\quad\quad\quad f_\eta^-(x) = \inf_{(g,g')\in B(1_G,\eta)^2} f(g^{-1}xg').$

Because of Inequality (16), there is a constant $c>0$ such that $R_{f_\eta^\pm}\leq c$, and for any $x$ we have $|f(x)-f_\eta^\pm(x)|\leq c\omega_f(c\eta)$. So there is (another) constant $c>0$ such that for any $\eta$ in $]0,1[$, we have

(19) $\quad\quad\quad\quad \int_{\mathrm{End}(V)} |f-f_\eta^\pm|d\mu_\infty \leq c\omega_f(c\eta).$

From the definition of $f_\eta^+$, and the fact that $\phi_\eta$ has support in $B(\Gamma,\eta)$, we can deduce that
$$\Psi(f\circ\iota_T)(\Gamma,\Gamma) \leq \int_{(\Gamma\backslash G)^2} \Psi(f_\eta^+\circ\iota_T)\alpha_\eta.$$
Thus, Lemma 3.2 applied to $f_\eta^+$ implies
$$\limsup_{T\to+\infty} \frac{\Psi(f\circ\iota_T)(\Gamma,\Gamma)}{T^d\ln(T)^e} \leq \int_{\mathrm{End}(V)} f_\eta^+ d\mu_\infty,$$

and Equation (19) implies

$$\limsup_{T\to+\infty} \frac{\Psi(f\circ\iota_T)(\Gamma,\Gamma)}{T^d \ln(T)^e} \leq \int_{\mathrm{End}(V)} f d\mu_\infty + c\omega_f(c\eta).$$

Since $\eta > 0$ is arbitrary, this shows the upper bound

$$\limsup_{T\to+\infty} \frac{\Psi(f\circ\iota_T)(\Gamma,\Gamma)}{T^d \ln(T)^e} \leq \int_{\mathrm{End}(V)} f d\mu_\infty.$$

Replacing $f_\eta^+$ by $f_\eta^-$, and reversing the inequalities, we obtain completly similar estimates and this concludes the proof of Theorem 2.

3.2. **Decay of Matrix Coefficients.** Here we establish the effective speed of convergence involved in the preceding proof in order to prove Proposition 2.2. For now on, we assume that moreover $G$ is almost simple of rank greater than 2. We will use H. Oh's estimates, described as follows. Let $\mathcal{S}$ be a *maximal strongly orthogonal system* (see [O]), the linear form

$$l = 1/2 \sum_{\alpha\in\mathcal{S}} \alpha \in \mathfrak{a}^+,$$

yields a family for $\theta > 0$,

$$\xi_\theta(kak') = c(\theta).\exp(-(1-\theta)\langle l|\log(a)\rangle),$$

where $g = kak'$ is the Cartan decomposition of a element $g$, and $c(\theta)$ a positive number sufficiently big. For all $\theta > 0$, any unitary representation $\pi$ of $G$ without invariant vectors, any two vectors $v, w$, Theorem 1.1 of [O] asserts that :

$$|\langle \pi(g)v, w\rangle| \leq \xi_\theta(g)(\dim\langle Kv\rangle)^{1/2}(\dim\langle Kw\rangle)^{1/2}||v||.||w||.$$

The following lemma is well-known.

**Lemma 3.3.** *Let $Q$ be a positive definite quadratic form on a real vector space $W$ of dimension $n$, and $\mathcal{Z}$ be a lattice in $W$. Then, provided $2k > n$, we have*

$$\sum_{v\in\mathcal{Z}-\{0\}} \frac{1}{Q(v)^k} < +\infty.$$

□

Let $\pi$ be a unitary representation of $G$ in a separable Hilbert space $\mathcal{H}$. Let

$$\mathcal{H} = \oplus_{\mu\in\hat{K}}\mathcal{H}_\mu,$$

be the decomposition into isotypical subrepresentations of $K$, that is $\mathcal{H}_\mu$ is a orthogonal sum of finitely or infinitely many copies of the same irreducible representation $\mu$ of $K$, which is finite dimensionnal.

Let $X_1,..,X_{dim(K)}$ be a orthonormal basis of $Lie(K)$ with respect to a $Ad$-invariant scalar product, and put $\Omega = 1 - \sum_{i=1}^{dim(K)} X_i^2$. This is a differential operator in the center of the envelloping algrebra of $Lie(K)$; thus it acts as multiplication by a scalar $c(\mu)$ on each isotypical component $\mathcal{H}_\mu$.

The following Theorem is due to Katok and Spatzier [KS], see also [KM]. The proof presented here is mainly a reproduction of theirs, but since we are interested in the effective constants involved, it seemed important to give all the details of the computations. Here $\Sigma_K^+$ denotes a set of positive roots for the Lie algebra of $K$.

**Theorem 6.** [KS] *Assume that $G$ is almost simple of rank greater than 2. There is a constant $C > 0$ such that for any $v, w$ $C^\infty$-vectors of a unitary representation without invariant vectors, and for any integer*

$$(20) \qquad 4m > \mathrm{rank}(K) + 2\#\Sigma_K^+,$$

*we have for all $g$ in $G$ and $\theta > 0$,*

$$|\langle \pi(g)v, w\rangle| \leq C.\xi_\theta(g)||\Omega^m v||.||\Omega^m w||,$$

*Proof.* From [Bo], §7, no 6, Proposition 4, there is a quadratic form $Q$, positive and definite on the dual of the Lie algebra of a Cartan subalgebra of Lie(K), such that if $\Lambda_\mu$ is a highest weight vector of $\mu$,

$$c(\mu) = 1 + Q(\Lambda_\mu + \rho) - Q(\rho),$$

where $\rho$ is the half sum of (a suitable choice of) positive roots $\Sigma_K^+$ of $Lie(K)$ (in fact this reference assumes that $Lie(K)$ is semisimple, but the formula still holds since $Lie(K)$ is a direct sum of an abelian and a semisimple Lie algebra). On the other hand, Hermann Weyl's formula asserts that

$$\dim(\mu) = \prod_{\alpha \in \Sigma_K^+} \frac{\langle \Lambda_\mu + \rho, \alpha \rangle}{\langle \rho, \alpha \rangle},$$

and is thus a polynomial of degree $\#\Sigma_K^+$ in $\Lambda_\mu$. This implies that there exists another quadratic form $Q_1$ such that

$$\dim(\mu) \leq Q_1(\Lambda_\mu)^{\#\Sigma_K^+/2}.$$

Thus there is a constant $c > 0$ such that for $\Lambda_\mu$ outside a small neighborhhood of zero, we have

$$c(\mu) \geq 1/cQ_1(\Lambda_\mu).$$

Let $v \in \mathcal{H}$ be a $C^\infty$ vector for $K$, that is a vector such that $g \mapsto \pi(g)v$ is a $C^\infty$ map, we can write

$$v = \sum_{\mu \in \hat{K}} v_\mu,$$

and for all $\mu$, $v_\mu \in \mathcal{H}_\mu$ is then also a $C^\infty$ vector. One have for all $m \geq 0$

$$||v_\mu|| = c(\mu)^{-m}||\Omega^m v_\mu||.$$

Let $v, w$ be two $C^\infty$ vectors, then

$$|\langle \pi(g)v, w\rangle| \leq \sum_{(\nu,\mu)\in \hat{K}^2} |\langle \pi(g)v_\mu, w_\nu\rangle|,$$

$$\leq \xi_\theta(g) \sum_{(\nu,\mu)\in \hat{K}^2} ||v_\mu||.||w_\nu||(\dim\langle Kv_\mu\rangle)^{1/2}(\dim\langle Kw_\nu\rangle)^{1/2}.$$

As is well-known, we have $\dim\langle Kv_\mu\rangle \leq \dim(\mu)^2$, and so we obtain

$$|\langle \pi(g)v, w\rangle| \leq \xi_\theta(g) \left(\sum_{\mu\in\hat{K}} ||v_\mu||\dim(\mu)\right) \left(\sum_{\nu\in\hat{K}} ||w_\nu||\dim(\nu)\right).$$

We have, for any integer $m \geq 1$,

$$|\langle \pi(g)v, w\rangle| \leq \xi_\theta(g) \left(\sum_{\mu\in\hat{K}} ||\Omega^m v_\mu||c(\mu)^{-m}\dim(\mu)\right) \left(\sum_{\nu\in\hat{K}} ||\Omega^m w_\nu||c(\nu)^{-m}\dim(\nu)\right),$$

and by Cauchy-Schwarz inequality,

$$|\langle \pi(g)v, w\rangle| \leq \xi_\theta(g)||\Omega^m v||.||\Omega^m v|| \left(\sum_{\mu\in\hat{K}} c(\mu)^{-2m}\dim(\mu)^2\right),$$

Now, since $c(\mu)^{-2m}\dim(\mu)^2$ is less than some constant times $Q_1(\Lambda_\mu)^{\#\Sigma_K^+ - 2m}$. Lemma 3.3 and the fact that highest weight vectors belongs to some lattice in the dual of the Lie algebra of a maximal torus in $K$ insures that

$$\sum_{\mu\in\hat{K}} c(\mu)^{-2m}\dim(\mu)^2 < +\infty,$$

provided that

$$4m > \text{rank}(K) + 2\#\Sigma_K^+.$$

□

### 3.3. Some test functions.

We fix $m$ to be equal to the least integer satisfying Equation (20). For $1 > \eta > 0$, we fix a map $\phi_\eta$ from $\Gamma\backslash G$ to $\mathbf{R}$, satisfying the following properties :

(1) $\phi_\eta$ is nonnegative and zero outside a ball of radius $\eta$ around $\Gamma$ in $\Gamma\backslash G$.

(2) $\phi_\eta$ is of class $C^\infty$, and $\int_{\Gamma\backslash G} \phi_\eta = 1$.

(3) The Sobolev norm of $W^{2m,2}(\Gamma\backslash G)$ of $\phi_\eta$ satisfies for some constant $C$ independant of $\eta$ between 0 and 1

$$||\phi_\eta||_{2m} \leq C\eta^{-2m-\dim(G)},$$

This is always possible (compare with Lemma 2.4.7 of [KM]). Since $\eta < 1$, we have also $||\phi_\eta - 1||_{2m} \leq (C+1)\eta^{-2m-\dim(G)}$. Now we put
$$\alpha_\eta(\Gamma g_1, \Gamma g_2) = \phi_\eta(\Gamma g_1)\phi_\eta(\Gamma g_2),$$

**Lemma 3.4.** *There exists a constants $C > 0$ such that for all $h$ in $G$, we have*
$$|\Psi^*\alpha_\eta(h) - 1| \leq C\xi_\theta(h)\eta^{-4m-2\dim(G)}.$$

*Proof.* As in the Lemma 3.1,
$$|\Psi^*\alpha_\eta(h) - 1| = |\langle \phi_\eta - 1, \pi(h^{-1})(\phi_\eta - 1)\rangle|,$$
and then apply Theorem 6, to obtain
$$|\Psi^*\alpha_\eta(h) - 1| \leq C\xi_\theta(h)||\Omega^m(\phi_\eta - 1)||^2.$$
Note that for some constant $C' > 0$, we have $||\Omega^m(\phi_\eta - 1)|| \leq C'||\phi_\eta - 1||_{2m}$. Together with the properties of $\phi_\eta$, this implies the desired estimate. □

### 3.4. Rate of weak convergence.

**Lemma 3.5.** *Let $u_\theta$ be the maximum on $\Delta$ of $\beta - (1-\theta)l$. Then for any $f$ in $\mathcal{B}$, we have if $T \geq 2/R_f$,*
$$\left|\int_{(\Gamma\backslash G)^2} \Psi(f \circ \iota_T)\alpha_\eta - \int_G f \circ \iota_T\right| \leq C||f||_\infty (R_f T)^{u_\theta} \ln(R_f T)^r \eta^{-4m-2\dim(G)}.$$

*Proof.* As in the proof of Lemma 3.2,
$$\left|\int_{(\Gamma\backslash G)^2} \Psi(f \circ \iota_T)(g, g')\alpha_\eta(g, g') dg dg' - \int_{\text{End}(V)} f \circ \iota_T d\mu_\infty\right|$$
$$\leq \int_G |(f \circ \iota_T)(h)(\Psi^*\alpha_\eta(h) - 1)| dh,$$
which is, by Lemma 3.4, less than
$$C\eta^{-4m-2\dim(G)} \int_G |(f \circ \iota_T)(h)|\xi_\theta(h) dh.$$
Lemma 2.15 applies to the integral and yields the expected result. □

### 3.5. Proof of Proposition 2.2.

**Lemma 3.6.** *Let $\lambda, \xi, \xi', \tau, h, \nu$ be given by Proposition 2.1. There is a constant $c > 0$ such that for any $\epsilon$ in $[0, 1[$, any $\eta$ in $]0, 1[$ and $f$ in $\mathcal{B}_1$, then $f_\eta^\pm$ has support in the ball of center 0 and radius $c$, and moreover we have for all $T \geq c\eta^{-\lambda}$,*
$$\left|\frac{1}{T^d \ln(T)^e} \int_G f_\eta^\pm \circ \iota_T d\mu - \int_{\text{End}(V)} f d\mu_\infty\right|$$
$$\leq c||f||_\infty \left(\nu(D_f(c\eta, \epsilon/2)) + \eta^{-\tau}T^{-\xi} + \frac{T^{-\xi'}}{\ln T} + h(\eta, T/c)\right) + c\epsilon.$$

*Proof.* From Inequality (16), there exists a $c > 0$ such that for $\eta < 1$ and all $g, g'$ in $B(1_G, r)$, we have
$$||g^{-1}xg' - x||_{\text{End}(V)} \leq c||x||_{\text{End}(V)}\eta,$$
and thus $f_\eta^\pm$ has support in the ball of radius $c+1$. So if $x$ in not in $D_f(c(c+1)\eta, \epsilon)$ and $||x||_{\text{End}(V)} \leq 1$, then $|f(x) - f_\eta^\pm(x)| \leq \epsilon$. So

$$\left|\int_G |f \circ \iota_T - f_\eta^\pm \circ \iota_T| d\mu_\infty\right| \leq 2||f||_\infty \int_{\text{End}(V)} 1_{D_f(c(c+1)\eta,\epsilon)} \circ \iota_T d\mu + \epsilon \int_{||x||_{\text{End}(V)} \leq T(c+1)} d\mu.$$

The last term can be bounded by some constant times $T^d \ln(T)^e$, thanks to Theorem 1. On the other hand, Proposition 2.1 applied to $f$ gives a bound of the required form for the difference $|\frac{1}{T^d \ln(T)^e} \int_G f \circ \iota_T d\mu - \int_{\text{End}(V)} f d\mu_\infty|$, so it is now sufficient to give a similar bound for $\int_{\text{End}(V)} 1_{D_f(c(c+1)\eta,\epsilon)} \circ \iota_T d\mu$. Let $g$ be the characteristic function of the set $D_f(c(c+1)\eta, \epsilon)$. It is a consequence of the triangle inequality that
$$D_g(\eta, 0) \subset D_f((c^2 + c + 1)\eta, \epsilon/2).$$
Thus, Theorem 2.1 applied to $g$ with $\delta' = \eta$ and $\epsilon' = 0$ implies that for some constant $C > 0$, we have for $T \geq C\eta^{-\lambda}$,
$$\frac{1}{T^d \ln T^e} \int_{\text{End}(V)} 1_{D_f(c(c+1)\eta,\epsilon)} \circ \iota_T d\mu \leq \mu_\infty(D_f(c(c+1)\eta, \epsilon))$$
$$+ C\left(\nu(D_f(C(c^2+c+1)\eta, \epsilon/2)) + \eta^{-\tau}T^{-\xi} + \frac{T^{-\xi'}}{\ln T} + h(\eta, T/C)\right).$$
Together with the facts that $\nu \geq \mu_\infty$ and $D_f(c(c+1)\eta, \epsilon) \subset D_f((c^2+c+1)\eta, \epsilon/2)$, this concludes the proof of the lemma. □

Let $f$ in $\mathcal{B}_1$. As before, we have
$$\int_{(\Gamma\backslash G)^2} \Psi(f_\eta^- \circ \iota_T)\alpha_\eta \leq \Psi f(\Gamma, \Gamma) \leq \int_{(\Gamma\backslash G)^2} \Psi(f_\eta^+ \circ \iota_T)\alpha_\eta.$$
Thus, Lemma 3.5 applied to $f_\eta^+$ implies that for some $c > 0$
$$\Psi f(\Gamma, \Gamma) \leq \int_G f_\eta^+ \circ \iota_T d\mu + c||f||_\infty (R_{f_\eta^+}T)^{u_\theta} \ln(R_{f_\eta^+}T)^r \eta^{-4m-2\dim(G)}.$$
Since $R_{f_\eta^+}$ is in fact bounded independently of $f$ in $\mathcal{B}_1$, we can be drop them in the expression, up to a change of the multiplicative constant $c > 0$. We can also bound $\ln(T)^r$ by some constant times $T^\theta$, provided $\theta > 0$. Thus, the use of Lemma 3.6 insures that for all $T \geq c\eta^{-\lambda}$, we have
$$\frac{\Psi f(\Gamma, \Gamma)}{T^d \ln(T)^e} - \int_G f d\mu_\infty \leq c\epsilon$$

$$+c||f||_\infty \left(\nu(D_f(c\eta,\epsilon/2)) + \eta^{-\tau}T^{-\xi} + \eta^{-4m-2\dim(G)}T^{u_\theta-d+\theta} + h(\eta,T/c)\right).$$

We now put $\tau_1 = 4m + 2\dim(G)$, $\xi_1 = d - u_\theta - \theta$ and $\delta = \eta$ to obtain the required form of Proposition 2.2.

The lower bound is similar. To conclude the proof of Proposition 2.2, one have to check that $\xi_1 > 0$; this is insured if $\theta$ is small enough by the fact that $u_\theta < d$ because $l$ is greater than half of any simple positive root, due to the fact that it is the half sum of a maximal orthogonal system.

## 4. Examples

**4.1. Standard representation of $SL(n,\mathbf{R})$.** Let us write for $\Phi = (\lambda_1,..,\lambda_n)$ the set of weights of the standard representation. Here $r = n - 1$. The roots are then $\Sigma = \{\lambda_i - \lambda_j\}_{i\neq j}$, and a set of positive roots is given by $\Sigma^+ = \{\lambda_i - \lambda_j\}_{i<j}$; the multiplicities $m_\alpha$ are in this case all equal to 1, and so, using $\sum \lambda_i = 0$, we have

$$\beta = \sum_{i=1}^{n-1}(2n-2i)\lambda_i.$$

It can be verified that in this case, $d = n(n-1)$, and $e = 0$, because $\beta/n(n-1)$ belongs to the face of $\mathcal{C}$ not containing $\lambda_n$, and to no other. Here the chosen basis $(\chi_1,..,\chi_{n-1})$ of $\mathfrak{a}^*$ is $(\lambda_1,..,\lambda_{n-1})$, and thus $\tau_i = d\mu_i - 1 = 2(n-i) - 1$. Moreover, the only $l_i$ not equal to zero is $l_n$, and $l_n = -n$, thus from Equation (11) we have $\lambda = 1/n$. Since $\lambda_n = -\chi_1 - ... - \chi_{n-1}$, so Equations (12) give

$$\tau = \sup_{1\leq i\leq n-1} -\frac{2(n-i)-1}{(n-1)(-1)} = \frac{2n-3}{n-1},$$

and

$$\xi = \inf_{1\leq i\leq n-1} \frac{(-n)(2(n-i)-1)}{(n-1)(-1)} = \frac{n}{n-1}.$$

Note $\Omega$ the subset such that

$$\prod_{\alpha\in\Sigma^+}(2\sinh\langle\alpha|a\rangle)^{m_\alpha} = \sum_{\gamma\in\Omega} h_\gamma \exp\langle\gamma|a\rangle.$$

Elements of $\Omega$ are obtained by chosing a sequence $\epsilon_{i,j} = \pm 1$ for $i < j$, every $\gamma \in \Omega$ being written

$$\gamma = \sum_{i<j}\epsilon_{i,j}(\lambda_i - \lambda_j) = \sum_i\left(\sum_{j>i}\epsilon_{i,j} - \sum_{j<i}\epsilon_{j,i}\right)\lambda_i.$$

The linear form $a$ defining $\mathcal{F}_\beta$ but the intersection of $\mathcal{C}$ with $\langle\chi|a\rangle = 1$ is $a = \sum_{i=1}^{n-1}\chi_i^*$. We compute

$$\langle\gamma|\sum_{i=1}^{n-1}\chi_i^*\rangle = n\left(\sum_{i<n}\epsilon_{i,n}\right).$$

Thus, since we know by Lemma 2.14 that $\gamma/(n^2 - n)$ is in $\mathcal{C}$, this determines the set $\Omega_1$ of $\gamma$ such that $\gamma/(n^2 - n)$ belongs to $\mathcal{F}_\beta$ : they are exactly the ones such that $\epsilon_{i,n} = 1$ for all $i = 1, .., n - 1$.

It remains to compute $\xi'$. The polyhedron $\mathcal{C}^*$ is the dual polyhedron of the simplex $\mathcal{C}$, whose extremal points are the point $x_i$ for $i = 1, .., n$ defined by $\lambda_j(x_i) = 1$ if $i \neq j$ and $-n + 1$ otherwise. Let $\gamma$ be in $\Omega - \Omega_1$, we compute

$$\langle \gamma + \beta | x_k \rangle / 2 = n \left( \sum_{i<k} \frac{\epsilon_{i,k} + 1}{2} - \sum_{j>k} \frac{\epsilon_{k,j} + 1}{2} \right).$$

Thus, we have

$$\langle \gamma + \beta | x_k \rangle / 2 \leq n \sum_{i<k} \frac{\epsilon_{i,k} + 1}{2}.$$

Remark that in all cases this is smaller than $n(n-2)$ because $\gamma$ does not belong to $\Omega_1$. Thus $(\gamma + \beta)/(2n(n-2))$ belongs to $\mathcal{C}$; and since $\gamma \leq (\gamma + \beta)/2$, the maximum of $\gamma$ on $\Delta$ is always less than $n(n-2)$. Thus we can take $\xi'$ any number satisfying $\xi' < n(n-1) - n(n-2)$, i.e.

$$\xi' < n.$$

Now we determine $\tau_1, \xi_1$.

The half sum $l$ of a maximal strongly orthogonal system in $A_{n-1}$ is given in Appendix of [O]. With our notations, if $n$ is even, we have

$$l = \frac{1}{2} \left( \sum_{i=1}^{n/2} \lambda_i - \sum_{i=n/2+1}^{n} \lambda_i \right),$$

and for $n$ odd we have

$$l = \frac{1}{2} \left( \sum_{i=1}^{(n-1)/2} \lambda_i - \sum_{i=(n-1)/2+2}^{n} \lambda_i \right).$$

The value of $\beta - l$ at an extremal points $x_k$ of $\mathcal{C}^*$ is

$$\langle \beta - l | x_k \rangle = n(2k - 1) - n^2 + \epsilon n/2,$$

with $\epsilon = 1$ if $k \leq n/2$, $\epsilon = -1$ if $k \geq n/2 + 1$ and $\epsilon = 0$ if $n$ is odd and $k = (n+1)/2$. Thus it is maximum when $k = n$, with maximum $n(n - 3/2)$; this proves that the maximum $u_0$ of $\beta - l$ on $\Delta$ is smaller than

$$u_0 \leq n(n - 3/2).$$

So we obtain that admissible $\xi_1$ are those satisfying $\xi_1 < d - n(n - 3/2)$, that is

$$\xi_1 < n/2.$$

A maximal compact subgroup in $SL(n, \mathbf{R})$ is $SO(n)$. We list in the following array for each values of $n$ the rank of $SO(n)$, the number of positive roots, a integer $m$ satisfying Equation (20) when Theorem 6 is applicable (i.e. $n > 2$) (we do not try to give an optimal one since it is of the same order than $\dim(G)$), and a real number that is greater or equal to $4m + 2\dim(G)$; to obtain the values indicated, see for example ([BD], V.6) for the second and third columns. In any case, $\dim(G) = n^2 - 1$.

| $n$ | rank$K$ | $\#\Sigma_K^+$ | $m$ | $\tau_1 \geq 4m + 2\dim(G)$ |
|---|---|---|---|---|
| 2 | 1 | 1 | $\emptyset$ | $\emptyset$ |
| 3 | 1 | 1 | 1 | 17 |
| 4 | 2 | 2 | 2 | 38 |
| 6 | 3 | 6 | 4 | 86 |
| $n \geq 5$ odd | $(n-1)/2$ | $(n^2 - 2n - 3)/4$ | $(n-1)^2/4$ | $3n^2 - 1$ |
| $n \geq 8$ even | $n/2$ | $(n^2 - 2n)/4$ | $n^2/4$ | $3n^2 - 1$ |

So we need now to minimize the upper bound given by Proposition 2.2. We introduce a small parameter $\theta > 0$, because the values $\xi' = n$ and $\xi_1 = n/2$ are a priori not allowed. We take $\epsilon = 0$ and $\delta = T^{-\xi_1/(\tau_1+1)}$, which is allowed since $1/\lambda = n \geq \xi_1/(\tau_1+1) = 1/(6n) + \theta$ if $n \geq 7$ and $n = 5$, and can be verified for $n = 3, 4, 6$. Thus the remainder term in Theorem 4 for the number of lattice points in a ball of radius $T$ is less than $O(T^{-\alpha})$, with $\alpha = \inf(\xi' = n/2 - \theta, 1/(6n) - \theta, (\xi(\tau_1+1) - \xi_1\tau)/(\tau_1+1))$, which is in this case $\alpha = 1/(6n) - \theta$ for $n \geq 7$. We summarize the results in the following proposition.

**Proposition 4.1.** *Let $\Gamma \subset SL(n, \mathbf{R})$ a lattice, with $n \geq 3$. If $n = 5$ or $n \geq 7$, we have for any $\theta > 0$*

$$\#\{\gamma \in \Gamma \mid ||\gamma|| \leq T\} = cT^{n^2-n} + O(T^{n^2-n-\frac{1}{6n}+\theta}),$$

*If $n = 3$, the remainder term is $O(T^{6-\frac{1}{12}+\theta})$, if $n = 4$ it is $O(T^{12-\frac{1}{20}+\theta})$, and if $n = 6$ it is $O(T^{30-\frac{1}{30}+\theta})$.*

As stated in the introduction, this remainder term is not as good as the one obtained if one uses the remainder terms obtained in [DRS] together with Theorem 3.

4.2. **Adjoint representation of $SL(3, \mathbf{R})$.** Let $\rho$ be the Adjoint representation of $G = SL(3, \mathbf{R})$. Let $E_{k,l} = (e_{i,j,k,l})_{i,j}$ be the $3 \times 3$ elementary matrices with $e_{i,j,k,l} = 1$ if $i = k$ and $l = j$, else 0. We fix the following basis of $V = \mathfrak{sl}(3, \mathbf{R})$:

$$(E_{1,3}, E_{1,1} - E_{3,3}, E_{2,2} - E_{3,3}, E_{1,2}, E_{2,3}, E_{2,1}, E_{2,3}, E_{3,2})$$

with respect to this basis, the weights of the action of the positve diagonal matrices $A$ of $SL(3, \mathbf{R})$ are

$$(\lambda_1 - \lambda_3, 0, 0, \lambda_1 - \lambda_2, \lambda_2 - \lambda_3, \lambda_2 - \lambda_1, \lambda_3 - \lambda_2, \lambda_3 - \lambda_2).$$

Here $d = 2$ and $e = 1$, because $\beta/2 = \lambda_1 - \lambda_3$ is a extremal point of the hexagon $\mathcal{C}$. It can be checked that up to a multiplicative constant, we have

$$\int_{\text{End}(V)} f d\mu_\infty = \int_{\mathbf{R}^+ \times K \times K} f\left(\rho(k)\text{diag}(x, 0, 0, 0, 0, 0, 0, 0)\rho(k')\right) x dx dk dk',$$

where the matrix representation is considered with respect to the given basis. Thus in this case, the support of the measure is a subset of the set of rank $\leq 1$ matrices (which is a closed set); however, the set of limit points of sequences $\rho(g_i)/T_i$ for $T_i$ tending to infinity contains rank 2 matrices, for example if we define $g_n$ by

$$\rho(g_n) = \rho(\mathrm{diag}(e^n, e^n, e^{-2n})) = \mathrm{diag}(e^{3n}, 1, 1, 1, e^{3n}, 1, e^{-3n}, e^{-3n}),$$

and, taking $T_n = e^{3n}$,

$$\lim_{n \to +\infty} \frac{\rho(g_n)}{e^{3n}} = \mathrm{diag}(1, 0, 0, 0, 1, 0, 0, 0),$$

which is of rank 2.

4.3. **Another example.** Here we consider $G = SL(3, \mathbf{R})$ acting on $V = \mathbf{R}^3 \oplus \mathbf{R}^3$ by $\rho(g)(v_1, v_2) = (gv_1, {}^t g^{-1} v_2)$. Keeping the notations of 4.1 and 4.2, the set of weights is

$$\Phi = \{\lambda_1, -\lambda_1, \lambda_2, -\lambda_2, \lambda_3, -\lambda_3\},$$

whose convex hull $\mathcal{C}$ is a regular hexagon. Here we fix the basis $\chi_1 = \lambda_1, \chi_2 = -\lambda_3 = \lambda_1 + \lambda_2$, thus $\beta = 4(\chi_1/2 + \chi_2/2)$, so $d = 4, e = 0, \tau_1 = \tau_2 = 1$. The other weights are $\chi_3 = \chi_1 - \chi_2, \chi_4 = -\chi_2, \chi_5 = -\chi_1, \chi_6 = \chi_2 - \chi_1$ with $l_3 = -1, l_4 = -2, l_5 = -2, l_6 = -1$ respectively, so we have by Equations (11) and (12)

$$\lambda = 1, \ \tau = 1, \ \xi = 1.$$

The set $\Omega$ contains 8 elements $\beta, \gamma_1 = 2\lambda_1 + 4\lambda_2, \gamma_2 = 2\lambda_1 - 2\lambda_2, 0, 0, -\gamma_1, -\gamma_2, -\beta$. The set $\Delta$ has extremal points $0, x_1, x_2, x_3$ with $\lambda_1(x_1) = \lambda_2(x_1) = 1/2, \lambda_1(x_2) = 1, \lambda_2(x_2) = -1/2$, and $\lambda_1(x_3) = 1, \lambda_2(x_3) = 0$. Computing the value of elements of $\Omega - \{\beta\}$ on these points gives a maximum value of 3, so we can take

$$\xi' < d - 3 = 1.$$

From the array of 4.1, we extract that $\tau_1 = 17$. Here $l = (\lambda_1 - \lambda_3)/2 = \beta/4$ so any $\xi_1$ satisfying $\xi_1 < 3$ is good.

Here is a interesting consequence of these calcultations.

**Proposition 4.2.** *Let $\Gamma \subset SL(3, \mathbf{R})$ be a lattice and $||.||$ be any norm on $\mathrm{End}(\mathbf{R}^3)$. There is some $c > 0$, depending on the norm, such that for any $\theta > 0$, we have*

$$\#\{\gamma \in \Gamma \mid ||\gamma|| \leq T, ||\gamma^{-1}|| \leq T\} = cT^4 + O(T^{4-1/6+\theta}).$$

4.4. **Counter-example.** Here we consider $\Gamma = SL(2, \mathbf{Z}) \times SL(2, \mathbf{Z}) \subset G = SL(2, \mathbf{R}) \times SL(2, \mathbf{R})$, and the representation $\rho$ of $G$ on $V = \mathbf{R}^2 \otimes \mathfrak{sl}(2, \mathbf{R})$. Let us call $G_1, G_2$ the first and second copy of $SL(2, \mathbf{R})$, and $\alpha_1, \alpha_2$ a nontrivial weight for the adjoint representation. Taking a basis $(e_1, e_2)$ of $\mathbf{R}^2$, the six weights of the representation $\rho_{|G_1}$ on $V$ are $\alpha_1/2, \alpha_1/2, \alpha_1/2, -\alpha_1/2, -\alpha_1/2, -\alpha_1/2$; the weights of $\rho$ are $\alpha_1/2 + \alpha_2, \alpha_1/2, \alpha_1/2 - \alpha_2, -\alpha_1/2 + \alpha_2, -\alpha_1/2, \alpha_1/2 - \alpha_2$, and here $\beta = \alpha_1 + \alpha_2$. Thus both $G$ and $G_1$ have growth rate $T^2$, and $G_2$ has growth rate $T$. We normalize Haar measure such that $SL(2,\mathbf{Z})$ has covolume 1. Let us write $\mu_{1,\infty}$ for the asymptotic limit of the Haar measure of $G_1$. Write $\Gamma_1 = SL(2,\mathbf{Z}) \times \{1\}$ and $\Gamma_2 = \{1\} \times SL(2,\mathbf{Z})$.

Then we have

**Proposition 4.3.** *Let $f$ be continuous of compact support in $\mathrm{End}(V)$. Then*

$$\lim_{T \to +\infty} \frac{1}{T^2} \sum_{\gamma \in \Gamma} f(\gamma/T) = \sum_{\gamma_2 \in \Gamma_2} \int_{\mathrm{End}(V)} f(\gamma_2 x) d\mu_{1,\infty}(x) < +\infty.$$

*Proof.* As usual, we assume that $f$ vanishes outside a ball of radius 1. Let $\gamma_2$ in $\Gamma_2$, we will first prove that for $T > 1$, we have for some constant $C > 0$ the following inequality

$$\left| \sum_{\gamma_1 \in \Gamma_1} f(\gamma_2 \gamma_1/T) \right| \leq C.T^2/||\gamma_2||^2.$$

We have $||\gamma_1 \gamma_2|| = ||\gamma_1|| \, ||\gamma_2||$ provided we consider the supremum norm of matrix coefficients in a base $(e_i \otimes f_j)_{i=1,2; f=1,2,3}$ of $V$, so the left-hand side is smaller than the number of $\gamma_1$ of norm less than $T/||\gamma_2||$, which is asymptotic to $c(T/||\gamma_2||)^2$ for some constant $c > 0$, by Corollary 1.2 applied to $\Gamma_1 \subset G_1$. We have

$$\frac{1}{T^2} \sum_{\gamma \in \Gamma} f(\gamma/T) = \sum_{\gamma_2 \in \Gamma_2} \left( \frac{1}{T^2} \sum_{\gamma_1 \in \Gamma_1} f(\gamma_2 \frac{\gamma_1}{T}) \right),$$

and each term of the sum over $\Gamma_2$ converge by Theorem 2 applied to $\Gamma_1 \subset G_1$ to

$$\int_{\mathrm{End}(V)} f(\gamma_2 x) d\mu_{1,\infty}(x).$$

We can apply the Lebesgue dominated convergence Theorem here because each term is dominated by $C/||\gamma_2||^2$ and

$$\sum_{\gamma_2 \in \Gamma_2} \frac{1}{||\gamma_2||^2} < +\infty,$$

this convergence being a consequence of the following asymptotic :

$$\#\{\gamma_2 \in \Gamma_2 : ||\gamma_2|| \leq T\} \sim c.T,$$

which is simply given by Corollary 1.2 applied to $\Gamma_2 \subset G_2$. □

We also have to show that the measure $\sum_{\gamma_2 \in \Gamma_2} (\gamma_2)_* \mu_{1,\infty}$ is different from the measure $\mu_\infty$ obtained by Theorem 1 applied to $G$. The measure $\mu_{1,\infty}$ has support on a 3-manifold and the measure $\mu_\infty$ is in the Lebesgue class of a manifold of dimension at least 4, thus cannot give positive measure to the support of $\mu_{1,\infty}$.

...

## 5. Acknowledgements

I would like to thank Antoine Chambert-Loir, Cornelia Drutu, Jean-François Quint, and Emmanuel Peyre for interest and useful comments.

École normale supérieure de Lyon, Unité de Mathématiques Pures et Appliquées, UMR CNRS 5669, 46, allée d'Italie, 69364 Lyon Cedex 07 France
*E-mail address*: `Francois.MAUCOURANT@umpa.ens-lyon.fr`